\documentstyle{amltd}
\begin{document}
\annalsline{151}{2000}
\received{March 16, 1998}
\startingpage{125}
\def\ritem#1{\item[{\rm #1}]}
\def\bye{\end{document}}
 
\font\tenrm=cmr10
\catcode`\@=11
\font\twelvemsb=msbm10 scaled 1100
\font\tenmsb=msbm10
\font\ninemsb=msbm10 scaled 800
\newfam\msbfam
\textfont\msbfam=\twelvemsb  \scriptfont\msbfam=\ninemsb
  \scriptscriptfont\msbfam=\ninemsb
\def\msb@{\hexnumber@\msbfam}
\def\Bbb{\relax\ifmmode\let\next\Bbb@\else
 \def\next{\errmessage{Use \string\Bbb\space only in math
mode}}\fi\next}
\def\Bbb@#1{{\Bbb@@{#1}}}
\def\Bbb@@#1{\fam\msbfam#1}
\catcode`\@=12

 \catcode`\@=11
\font\twelveeuf=eufm10 scaled 1100
\font\teneuf=eufm10
\font\nineeuf=eufm7 scaled 1100
\newfam\euffam
\textfont\euffam=\twelveeuf  \scriptfont\euffam=\teneuf
  \scriptscriptfont\euffam=\nineeuf
\def\euf@{\hexnumber@\euffam}
\def\frak{\relax\ifmmode\let\next\frak@\else
 \def\next{\errmessage{Use \string\frak\space only in math
mode}}\fi\next}
\def\frak@#1{{\frak@@{#1}}}
\def\frak@@#1{\fam\euffam#1}
\catcode`\@=12

  \input amssym.def
\input amssym.tex

\title{Linking numbers and boundaries\\ of varieties}
 \shorttitle{Linking numbers and boundaries of varieties} 

\twoauthors{H. Alexander}{John Wermer}

\institutions{University of Illinois at Chicago, Chicago, IL\\
\vskip6pt
Brown Universtiy, Providence, RI\\
{\eightpoint {\it E-mail address\/}: wermer@math.brown.edu}}

\newcommand{\LL}{{ \kern .25ex \rule{0.10ex}{1.5ex}%
 \rule{1.3ex}{0.10ex}\kern .25ex }}

\intro
  The intersection
 index  at a common point
 of two analytic varieties of complementary dimensions
 in
 $\Bbb C^n$  is positive.  This  observation, which   has been  called a
 ``cornerstone'' of algebraic geometry ([GH, p.~62]),  is a simple
 consequence of the fact that analytic varieties carry a natural
 orientation.
 Recast in terms of linking numbers,
   it is our  principal motivation.  It implies the following:
    Let $M$ be a smooth oriented compact
   3-manifold in $\Bbb C^3$.  Suppose  that $M$ bounds a
   bounded complex 2-variety  $V$.
    Here ``bounds''  means, in the sense of Stokes'
   theorem, i.e., that
   ${ b[V]}={ [M]}$ as currents.  Let $A$ be an algebraic curve in $\Bbb C^3$
   which is disjoint from M.
   Consider   the linking number  ${\rm link}(M,A)$ of
   $M$ and $A$.
   Since this  linking number is equal to
    the intersection number (i.e.\ the sum of the
   intersection indices) of
   $V$ and $A$,  by the
   positivity of these intersection  indices,
   we have  ${\rm link}(M,A) \geq 0$. The linking number
   will of course be
   $0$ if  $V$ and $A$
   are disjoint.  (As $A$ is not compact, this usage  of
   ``linking number''  will be clarified later.) This reasoning
   shows more generally that ${\rm link}(M,A) \geq 0$ if
   $M$ bounds a positive  holomorphic  2-chain.
   Recall that a  {\it  holomorphic $k$-chain\/}  in
   $\Omega \subseteq \Bbb C^n$ is a sum $\sum n_j [V_j]$ where
   $\{V_j\}$ is a locally finite family of irreducible $k$-dimensional
   subvarieties of $\Omega$ and $n_j \in \Bbb Z$ and that the
   holomorphic 2-chain  is
   {\it positive\/} if  $n_j >0$ for all $j$.
   Our first  result is that, conversely, the nonnegativity
   of the linking  number
   characterizes
   boundaries of
    positive holomorphic $2$-chains. 

\nonumproclaim{Theorem 1} 
Let $M$ be a smooth{\rm ,} oriented{\rm ,} compact{\rm ,}
   $3$\/{\rm -}\/manifold {\rm (}\/not necessarily connected\/{\rm )}  in
   $\Bbb C^3${\rm .} Suppose that ${\rm link}(M,A) \geq 0$
    whenever $A$   is  an algebraic curve in $\Bbb C^3$
   disjoint from $M${\rm .}  Then  there exists a  {\rm (}\/unique\/{\rm )}    positive
    holomorphic $2$\/{\rm -}\/chain $T$
   in $\Bbb C^3 \setminus M$ of finite mass
   and with bounded support  such that
    $[M]=b[T]${\rm .  }
\endproclaim

   We shall refer to the linking hypothesis in Theorem 1 as the
linking
   condition. More generally, a  smooth oriented compact manifold
   $M$ in $\Bbb C^n$ of (odd)
real  dimension  $k$  satisfies the {\it  linking condition\/} if
 $${\rm link}(M,A) \geq 0$$
 for  all algebraic subvarieties  $A$ of $\Bbb C^n$
disjoint from $M$ of pure (complex) dimension
$n -(k+1)/2$.
    Of course, the conclusion of Theorem 1
     is   closely related to the fundamental  result of Harvey
   and Lawson [HL] that  $M$ bounds a
   bounded holomorphic  2-chain $T$
   if and only if $M$ is maximally complex.
   In Theorem 1---unlike   the    Harvey-Lawson theorem---the
     holomorphic $2$-chain is positive.
   This reflects the fact that ``maximal complexity of $M$'' is
   unaffected by a change of
   orientation of $M$,  while our hypothesis  on
   linking numbers
    is tied to a specific orientation.
   One of the  main steps in our proof of
   Theorem 1 is indeed to verify that
   $M$  is maximally complex.

   If  $M$ bounds a holomorphic 2-chain $T=\sum n_j[V_j]$
    as in the last theorem, then
   by the maximum   principle
   $ \hbox{\rm supp }T \setminus M
   = \bigcup  V_j  \subseteq  \hat{M}$. Here
   $\hat{K} $, the polynomially convex hull
   of a compact set  $K \subseteq  \Bbb C^3$,  is defined
   as the set $\{z \in \Bbb C^3 :  |P(z)|  \leq   \sup_K  |P|
   \mbox{ for all polynomials } P  \mbox{ in } \Bbb C^3\}$.
   In general,  $\hat{M}$ will be larger that $\bigcup  V_j  \cup M$.
   While the points in the polynomial hull  are given explicitly by the
   definition  just stated,
   the  ``individual'' points of $ \hbox{\rm supp } T$,  with $T$
   the (unique) Harvey-Lawson solution to the
   equation $b[T]=[M]$
   for a given maximally complex $M$, on the other hand,
   are not explicitly given. The next result determines
   these points in terms of linking numbers. 
\vglue6pt

   \nonumproclaim{Theorem 2}  Let $M$ be as given  in Theorem {\rm 1} and let
   $T$ be the unique  bounded holomorphic
   {\rm 2-}\/chain in $\Bbb C^3$  such that  $b[T]=[M]${\rm .}
   Then for $ x \in \Bbb C^3 \setminus M${\rm ,}  $x \in  \hbox{\rm supp }T$
    if and only if

 $${\rm link}(M,A) >0$$
  \vglue6pt
\noindent  for every algebraic curve $A$ in $\Bbb C^3$ such
   that $x \in A$ and $A  \cap M = \emptyset${\rm .}
\endproclaim
\vglue6pt

   Of course, half of this equivalence is trivial: it is merely the
   above-mentioned positivity
   of the intersection numbers.  For the opposite implication,
   we shall show that if   $x \not\in  \hbox{\rm supp }T$,
    then there exists $A$ such that
   $x \in A$, $A  \cap M = \emptyset$, and
   ${\rm link}(M,A) =0$. 

  In order to prove these two theorems  about $3$-manifolds in $\Bbb
C^3$,
   we need to establish the
  corresponding theorems for smooth oriented 1-manifolds
  $\gamma$, that  are compact, but not
  necessarily connected.  Thus $\gamma$ is a
  finite disjoint union of  oriented   simple closed curves in $\Bbb
C^n$.
  Recall that
  $\gamma$ satisfies the {\it moment condition \/}  if  
  $$ \int_\gamma   \phi =0$$\vglue6pt
\noindent 
  for all holomorphic $(1,0)$-forms  $\phi$  in $\Bbb C^n$.
  By Harvey and Lawson [HL], if
  $\gamma$ satisfies the   moment condition, then
  $\gamma$ bounds  (in the sense of Stokes' theorem) a
  (unique) bounded  holomorphic
   $1$-chain in $\Bbb C^n \setminus \gamma$. 

  \nonumproclaim{Theorem 3}
  Let $\gamma$ be a smooth
  compact  oriented $1$\/{\rm -}\/manifold in $\Bbb C^n${\rm .}
  Suppose that ${\rm link}(\gamma,A) \geq 0$ for every
  algebraic hypersurface $A$ in $\Bbb C^n$ such that
  $A  \cap \gamma = \emptyset${\rm .}  Then $\gamma$ satisfies the
moment
  condition and  there exists a
  {\rm (}\/unique\/{\rm )}   positive holomorphic
   $1$\/{\rm -}\/chain $T$ in $\Bbb C^n \setminus \gamma$
   of bounded support and finite mass
   such that
  $b[T]=  [\gamma]$.{\rm }
\endproclaim

  Theorem 3 is the direct analogue for curves of Theorem 1. The
analogue
  for Theorem 2 is the following:

  \nonumproclaim{Theorem 4}
Let $\gamma$ be given as in Theorem {\rm 3} and
let
   $T$ be the unique  bounded positive homomorphic $1$\/{\rm -}\/chain
   in $\Bbb C^n \setminus \gamma$ such that  $b[T]=[\gamma]${\rm .}
   Then for $ x \in \Bbb C^n \setminus \gamma${\rm ,}
   $x \in  \hbox{\rm supp }T$ if and only if

  $${\rm link}(\gamma,A) >0$$
 \vglue6pt\noindent  for every algebraic hypersurface $A$ in $\Bbb C^n$ such
   that $x \in A$ and $A  \cap \gamma = \emptyset${\rm .}
\endproclaim

    Theorems 1 and 2 of course suggest that  corresponding
   results might   hold  for manifolds of odd dimension
    in all  $\Bbb C^n$. This turns out to be true.
   However, the main  steps
   in proving the general result are to establish the
   preliminary cases stated so
   far. For this reason, we have stated them
   separately, even though they are special cases
   of the general result, which we now state. 

    \nonumproclaim{Theorem 5}
 Let $M$ be a smooth oriented compact
manifold
in $\Bbb C^n$ of {\rm (}\/odd\/{\rm )}
real  dimension  $k$ with $3 \leq k \leq  2n-3${\rm .}  Then  $M$ satisfies
the
linking condition\/{\rm :}\/
$${\rm link}(M,A) \geq 0$$
\vglue6pt\noindent
for all algebraic subvarieties  $A$ of $\Bbb C^n$
disjoint from $M$ of pure {\rm (}\/complex\/{\rm )} dimension
$n -(k+1)/2$ if and only if  $M$ is maximally
complex and   there exists a {\rm (}\/unique\/{\rm )} positive
holomorphic $k$\/{\rm -}\/chain
 $T$
of dimension $ (k+1)/2 $ in $\Bbb C^n \setminus M$
of finite mass and bounded support such that  $[M]=b[T]${\rm .}
Moreover{\rm ,}  for all $x \in \Bbb C^n \setminus M${\rm ,}
$x \in  \hbox{\rm supp }T$ if and only if 

$${\rm link}(M,A)  > 0$$
\vglue6pt\noindent
for all algebraic subvarieties  $A$ of $\Bbb C^n$
disjoint from $M$ of pure {\rm (}\/complex\/{\rm )} dimension
$n -(k+1)/2$ such that $x \in A${\rm .}
\endproclaim

It may be of interest to reformulate two of our results. 

\medbreak
I) (Theorem 3 +Lemma 1.2). Let $\gamma$ be a smooth
oriented compact curve
in $\Bbb C^n$. Then there exists a positive
holomorphic $1$-chain $V$ in
$\Bbb C^n \setminus \gamma$ of finite mass such that
$[\gamma]=b[V]$  if and only if 
$\frac{1}{2\pi i} \int_\gamma \frac{dp}{p } \geq 0$ for any
polynomial
$p$ in $\Bbb C^n$ such that $p|_\gamma  \not = 0$.\medbreak

The second much simpler part consists of the
following statement, which
is proved, but not explicitly formulated, below.
It was conjectured by Dolbault and
Henkin ([DH, p.~388]) and was recently also proved by Dinh [D].
(We thank the referee for these references.)\medbreak

II) Let $M$ be a smooth compact manifold in $\Bbb C^n$
of real dimension
$2p-1\break\geq~3$. Then there exists a holomorphic $p$-chain
$V$ in $\Bbb C^n
\setminus M$ of finite mass such that
$[M]=b[V]$ if and only if  for almost any complex $(n-p+1)$-plane
$H$ in $\Bbb C^n$
the curve $\gamma = H \cap M$ bounds a holomorphic $1$-chain in
$H \setminus\gamma$.\medbreak

We shall begin with some preliminary remarks; these
include material  on linking numbers,
polynomial 
hulls of curves and the Arens-Royden theorem.
We then  establish first
the theorems for
curves. This relies on the theory of polynomial hulls of curves due to
Wermer [W], Bishop and Stolzenberg [St]as well as on
the Harvey-Lawson [HL] theorem for curves  that
involves the moment condition.
This is the most difficult case, at least in the
smooth case, in part  because
we do not know that $\hat{\gamma}$ is a  `nice' topological space for
the most  general smooth
$\gamma$.  When $\gamma$ is real analytic, then
$\hat{\gamma}$, as a topological space, is a finite simplicial complex,
and the proof is much
shorter than for  the smooth case.
From the curve result we deduce the theorem for
 $3$-manifolds in $\Bbb C^3$.  The
remaining cases are  then  obtained
by using projections  for $k=3$ and,
for  higher dimensions,   by slicing and an inductive procedure.
We shall use some standard facts and the notation  for currents; for this
we
refer to Federer [F], Harvey[H] and Harvey-Shiffman[HS].
To avoid confusion with
other uses of $\partial$, we denote
the boundary of a current $T$ by $bT$.
Hausdorff $k$-dimensional measure will be denoted by
${\cal H}^k$. We want to thank
Bruno Harris for some helpful conversations
on algebraic topology. 

\section{Preliminaries}

\demo{{\rm A.} Linking}
 We shall briefly recall the definition of linking
number
and then derive a few of its properties. For more details
we refer to
Bott and Tu  [BT] who give  enlightening
discussions of  the linking number ([BT, pp.~231--235])
and the Poincar\'{e}  dual.
A very  general definition of linking number for
singular homology classes is  given by
Spanier  ([Sp, p.~361]).
Let $M$ and $Y$ be
disjoint compact smooth oriented submanifolds of
$\Bbb R^N$ of respective dimensions  $s$ and $t$. Suppose that
$s+t=N-1$.
Then the linking number ${\rm link}(M,Y)$ can be defined as
follows: Let  $\Sigma$ be a compact  oriented
$(s+1)$-chain in $\Bbb R^N$ such that $M=b\Sigma$
and such that $\Sigma$ and $Y$ meet transversally.
Then  ${\rm link}(M,Y)$ is the intersection number
$\#(\Sigma,Y)$.

There is  a useful alternate equivalent definition
of the linking number which uses the Poincar\'{e}  dual:
Let $\eta_M$  and $\eta_Y$ be compact Poincar\'{e}  duals
of $M$ and $Y$ supported on disjoint neighborhoods of $M$
and $Y$ respectively. Thus  $\eta_M$  is a closed $(N-s)$-form with
compact
support in $\Bbb R^N$ and so there exists a compactly supported
$(N-s-1)$-form $\omega_M$ in $\Bbb R^N$ such that
$d \omega_M=  \eta_M $. Then
$$ {\rm link}(M,Y)= \int \omega_M \wedge  \eta_Y. \leqno(1.1)$$
The integration on the right-hand side of (1.1) is over all of 
$\Bbb R^N$, but of course  $\omega_M \wedge  \eta_Y$ has
 compact support.   The Poincar\'{e} dual
 $\eta_Y$ can be ``localized'' to have
 support in an arbitrary neighborhood $W$ of $Y$. Its fundamental
property is
 that
 $\int_W  \phi \wedge \eta_Y  = \int_Y \phi$
 for all {\it  closed\/}  $t$-forms $\phi$   on $W$.
 Choosing  $\phi$ to be the restriction of $\omega_M$ to $W $
($\omega_M$ is
 closed as a form on $W$)  we get
 $$ {\rm link}(M,Y)= \int \omega_M \wedge  \eta_Y=
 \int_W \omega_M \wedge  \eta_Y=\int_Y \omega_M.\leqno(1.2)$$

    We shall use the linking number in a somewhat more general
setting.
    Namely,  we need ${\rm link}(M,A)$ when $M$ is a compact
oriented
    $k$-manifold ($k $ odd) in $\Bbb C^n$ and $A$ is an algebraic
subvariety of
    $\Bbb C^n$ with its natural orientation
    and of complex dimension $s$ so that $k+2s=2n-1$.
    Then $A$ is not compact and the above definitions of linking
number
    need to be extended.
    One approach is to  modify $A$ outside of a large
    ball  $B(r)$, centered at $0$ of radius $r$,  and
    containing $M$,  so that $A$  becomes compact  as
    follows:  Let $R = A \cap bB(r)$, a compact oriented $(2s-1) $-
chain 
    (for almost all $r$) contained in the sphere $bB(r)$,
    and let $A''$ be a $2s$-chain in
    $bB(r)$ so that $bA'' = R$. Then
    $A'= A \cap B(r)   -  A''$ is a  (compact)
    $2s$-cycle in $\Bbb C^n$ which agrees with $A$ inside $B(r)$.
    We  take  ${\rm link}(M,A)$  to be ${\rm link}(M,A')$; it
    is independent of the choices of $r$ and $R$.
    Alternatively,  we   can apply the
    first definition above,
    taking   ${\rm link}(M,A)$  as $\#(\Sigma,A)$ where
    $M=b\Sigma$
and   $\Sigma \subseteq B(r)$  is such that $\Sigma$ and
 $A$ meet
 transversally.  This yields the same linking number.
 The definition in terms of differential forms can also be adapted to
this
 setting as follows. Let $[A]$ be the current of integration over $A$,
 a positive $(s,s)$-current. We can extend the  definition of (1.2) to the
 following:
 $$ {\rm link}(M,A)=   \int_A\omega_M=
[A](\omega_M).\leqno(1.3)$$ 
\enddemo

\nonumproclaim{Lemma 1.1} Let $M$ be a smooth
real $k$ dimensional compact oriented
manifold in $\Bbb C^n$ and let $H$ be a complex
hyperplane in $\Bbb C^n$
given as $\{F=\lambda\}${\rm ,} where $F$ is a
complex linear function on $\Bbb
C^n${\rm ;} we view $H$ as a copy of $\Bbb C^{n-1}${\rm .}
  Suppose that $Q= M \cap H$ is a smooth $k-2$ manifold{\rm ,}
  oriented as the
slice of $M$ by  the map $F${\rm .}  Let $A$ be an algebraic variety  of
pure  complex  dimension
$n-(k+1)/2$ contained in $H$ and disjoint from $Q${\rm .}
Then
${\rm link}(M,A)${\rm ,} the {\rm ``}\/link\/{\rm ''} taken in $\Bbb C^n${\rm ,} agrees with
${\rm link} (Q,A)${\rm ,} the {\rm ``}\/link\/{\rm ''} taken in $H$ and well\/{\rm -}\/defined in $H$
since
$2n-k-1= 2(n-1) - (k-2) -1${\rm .}
\endproclaim

{\it Proofs}.  (1) Let $j: H \rightarrow \Bbb C^n$ be the inclusion map. Let
$\eta_M$ be a  Poincar\'{e} dual of $M$ in $\Bbb C^n$.
We can choose $\eta_M$ to have compact support disjoint from $A$.
By a basic
functorial property of Poincar\'{e} duals ([BT, p.~69]) we have
$j^*(\eta_M)=\eta_{j^{-1} (M)} = \eta_{M \cap H} =\eta_Q$,
a   Poincar\'{e} dual of $Q$ in $H$.
Let $\omega_M$ be a compactly supported
 $(2n-k-1)$-form in $\Bbb C^n$ such that $d\omega_M=\eta_M$.
 Set $\omega_Q= j^*(\omega_M)$. Now, $d\omega_Q=d(j^*(\omega_M))=j^*(d\omega_M)=
 j^*( \eta_M)= \eta_Q$.  Hence, by two applications of (1.3),
\medbreak
\centerline{\hfill${\displaystyle {\rm link}(Q,A)= \int_{j^{-1}(A)}  \omega_Q =
\int_{j^{-1}(A)}  j^*(\omega_M)=
\int_A \omega_M= {\rm link}(M,A).}$ \hfill\qed}
\bigbreak

(2) Let $G$ be a $(k+1)$-chain in $\Bbb C^n$ such that $bG=M$ and such
that $G$   and $G\cap H$ intersect $A$ transversally.
One checks that $\#(G,A)$ in $\Bbb C^n$ equals
$\#(G\cap H,A)$ in $H$.  This implies that
${\rm link}(Q,A)={\rm link}(M,A)$. \hfill\qed
\bigbreak

Let $\gamma$ be a smooth $1$-cycle in $\Bbb C^n$ and let
$A$ be an algebraic hypersurface in $\Bbb C^n$ that is disjoint from
$\gamma$.

\nonumproclaim{Lemma 1.2}  If $A=Z(P)${\rm ,} where $P$ is a polynomial in
$\Bbb C^n${\rm ,} then
$${\rm link}(\gamma,A)= \frac{1}{2\pi i} \int_\gamma  dP/P
.$$
\endproclaim

{\it Proofs}.  (1) We give first a proof based on the Poincar\'{e}-Lelong
formula
$$[A] = -i/2\pi d\partial \log |P|^2,$$  where $[A]$ is the
$(n-1,n-1)$-current
of integration over $A$.
Then $\psi=-i/2\pi  \partial \log |P|^2  $ is a current
such that $d \psi= [A]$.
Off of the  zero set $A=Z(P)$, we have $\psi=\frac{1}{2\pi i}  dP/P.$
 We can obtain a smooth form cohomologous to
   $[A]$ ([GH, p.~393]), by taking the convolution of
   $[A]$ with a smooth function
 and this smooth form then is a Poincar\'{e} dual $\eta_A$ to $A$.
 Corresponding to $\eta_A$ is a smoothing  $\omega_A$ of
 $  \psi$ such
 that $d\omega_A= \eta_A$ and such that $\omega_A$ is
 cohomologous to    $  \psi$.
 Let $\eta_\gamma$ be a compact Poincar\'{e} dual of $\gamma$,
 a $(2n-2)$-form, supported
 on a small neighborhood of $\gamma$  that  is disjoint
 from the support of $ \eta_A$. We have, since
 $\omega_A$  is
 cohomologous to  $\frac{1}{2\pi i}\frac{dP}{P}$ off of $A$,
 
\medbreak
\centerline{\hfill ${\displaystyle {\rm link}(\gamma,A)=
 {\rm link}(A,\gamma)=
 \int \omega_A    \wedge  \eta_\gamma=
 \int_\gamma \omega_A=
 \frac{1}{2\pi i} \int_\gamma  \frac{dP}{P}.}$\hfill\qed}
\bigbreak

  (2) Consider the map $\psi: \Bbb C^n \rightarrow \Bbb C^{n+1}$ given
by
$\psi(z)=(z,P(z))$. Set $\gamma'=\psi(\gamma)$ and
$A'= \{w \in \Bbb C^{n+1}: w_{n+1}=0\}$.
Then ${\rm link}(\gamma,A)$ in $\Bbb C^n$ equals
${\rm link}(\gamma', A')$ in $\Bbb C^{n+1}$. One can
continuously deform $\gamma'$ in  $\Bbb C^{n+1}$ to
the curve $\gamma''= (0 \in\Bbb C^n)  \times P(\gamma)$
in  $\Bbb C^{n+1}$ by
curves $\gamma_t= \{(tz,P(z)): z \in \gamma \}$, $0 \leq t \leq1$,
that are disjoint from $A'$ in  $\Bbb C^{n+1}$.
Hence  ${\rm link}(\gamma', A')=
{\rm link}(\gamma'',A')$ in $\Bbb C^{n+1}$.
Finally ${\rm link}(\gamma'',A')$  in $\Bbb C^{n+1}$
equals ${\rm link}(P(\gamma),\{0\})$
in $\Bbb C$ and this last linking number in $\Bbb C$ is just
the winding number of $P(\gamma)$ about $0$ which is
$\frac{1}{2\pi i} \int_\gamma  \frac{dP}{P}$. \hfill\qed\medbreak

\demo{{R}emark}   It may be of interest to observe, although we shall
not need it, that
Lemma 1.2 extends to the higher dimensional setting of
Theorem 5.  Namely, suppose that $M$ has dimension $k$
and that $A$, of complex dimension $s$,
where $2s +k =2n-1$,
 is a complete intersection in $\Bbb C^n$
given as the common zero set of polynomials $P_1,P_2, \cdots, P_{n-
s}$.
Let
$P= (P_1,P_2, \cdots, P_{n-s}): \Bbb C^n \rightarrow  \Bbb C^{n-s}$
and let $\beta_{n-s}$ be the Bochner-Martinelli [GH]
$(2(n-s)-1)$-form in $\Bbb C^{n-s}$
with singularity  at $0$. Then $k=2(n-s)-1$ and
$${\rm link}(M,A)= \  \int_M  P^*(\beta_{n-s}).$$
This can be easily verified by adapting the
second proof of Lemma 1.2.\enddemo

\demo{{\rm B.} Polynomial hulls of curves} For $K$ a compact subset
of $\Bbb C^n$,  ${\bf P}(K)$  will denote the uniform closure on $K$
of the polynomials  in $\Bbb C^n$
 and  $\hat{K}$ will denote the polynomially convex hull of $K$,
defined as the set $\{z \in \Bbb C^n:  |f(z)| \leq \sup_K |f| \mbox{ for
all
polynomials }  f  \mbox{ in }   \Bbb C^n \}$. The maximal ideal space
of ${\bf P}(K)$ can be identified with $\hat{K}$. Then the Shilov
boundary of
${\bf P}(K)$ is identified with a subset of $K$.\enddemo

\nonumproclaim{Lemma  1.3} Let $\Gamma_1$ be a finite union of smooth
curves
 in $\Bbb C^n$ and let $\beta$ be a smooth arc in $\Bbb C^n$ which is
 disjoint from $\Gamma_1${\rm .} Then
 $$\widehat{(\Gamma_1 \cup \beta})= \widehat{\Gamma_1} \cup
\beta.$$
\endproclaim

 \demo{Proof} We need only to show that
 $\widehat{(\Gamma_1 \cup \beta})
 \subseteq \widehat{\Gamma_1} \cup \beta$,
 the opposite inclusion being trivial.
 We know by Stolzenberg [St] that
 $ \widehat{ \Gamma_1
 \cup \beta} \setminus  \Gamma_1 \cup \beta$
 is a possibly empty, one-dimensional subvariety  $V$ of
 $\Bbb C^n \setminus  \Gamma_1 \cup \beta$.
 We claim that $V \subseteq \widehat{\Gamma_1}$.
 Suppose not.  Then there exists a polynomial
 $f$ in $\Bbb C^n$ such that $|f| < 1/2$ on $\widehat{\Gamma_1}$
 and  $f(p)=1$ for some $p \in V$. We can adjust $f$ so that
 $f \not = 1$ on
 $\beta$.  Set $g=1-f$. Then $\mbox{Re}(g) >0$ on
 $ \Gamma_1 $ and so $g$ has a continuous logarithm on
 $ \Gamma_1$. As $g\not = 0$ on $\beta$, $g$ also
 has a continuous logarithm  on the arc $\beta$. Hence,
 $\Gamma_1$ and $\beta$ being disjoint,
 $g$ has a logarithm on
 $\Gamma_1 \cup \beta \supseteq bV=\bar{V} \setminus V   $.
 Thus, by the argument
 principle [St], $g$ has no zeros on $V$. But $g(p)=0$, a contradiction.
 Hence  the claim $V \subseteq \widehat{\Gamma_1}$.
 Therefore $\widehat{\Gamma_1 \cup \beta} \subseteq
 \widehat{\Gamma_1} \cup \beta$ and the lemma
 follows.\enddemo

\nonumproclaim{Lemma 1.4}  Let $\Gamma$ be a finite union of
smooth disjoint simple
closed curves in $\Bbb C^n${\rm .} Suppose that
\begin{itemize}
\ritem{(a)}  $\Gamma$ is contained in the closure of
$\hat{\Gamma} \setminus \Gamma${\rm ,} and 

\ritem{(b)}  $\Gamma$ is the Shilov boundary of ${\bf P}(\hat{\Gamma})${\rm .}
\end{itemize}
\noindent 
Let $E$ be the complement in $\Gamma$ of
the set of points  $p \in \Gamma$ such that
the pair $(\hat{\Gamma},\Gamma)$ is locally a smooth $2$\/{\rm -}\/manifold
with boundary  {\rm (}\/contained in $\Gamma$\/{\rm )} in a neighborhood of
$p${\rm .}
Then  $E  \subseteq \Gamma$ is compact with
${\cal H}^1(E)=0${\rm .}
\endproclaim

{\it Proof}. The set  of points of $\Gamma$
where $(\hat{\Gamma},\Gamma)$ is
locally a smooth $2$-manifold
is open in $\Gamma$  and so $E$ is compact.
  By [St], $\hat{\Gamma} \setminus \Gamma$
is  a nonempty $1$-dimensional
subvariety of $\Bbb C^n \setminus \Gamma$.
We argue by contradiction and suppose that ${\cal H}^1(E)>0$.
Let $p \in E$
 be a point of density in $\Gamma$ of $E$. Choose a polynomial
 $f$ such that  $p$ is a regular point of $f|\Gamma$.
 Hence there is a subarc $\tau$ of $\Gamma$  such
 that ${\cal H}^1(\tau \cap E)>0$
  and  such that $f$ maps $\tau$ diffeomorphically to an arc $\tau'
 \subseteq \Bbb C$.  Since the set of singular values of $f|\Gamma$
has
 ${\cal H}^1$-measure zero, by shrinking
  $\tau$ and $\tau'$ we can further
 assume that $\tau'$ contains no singular values
  of $f|\Gamma$ and
 that $f^{-1}( \tau') \cap \Gamma$ is the union
 of $s$ arcs $\tau_1,\tau_2, \cdots \tau_s$ such
 that $\tau_1 =\tau$ and each $\tau_j$ is
  mapped by $f$  diffeomorphically to   $\tau' $.
  Choose a small neighborhood $\omega $ of $\tau'$ in
  $\Bbb C$ such that
  $f(\Gamma) \cap \omega = \tau'$ and $\omega \setminus \tau'$
  is the union of two
  components $\Omega_1$ and $\Omega_2$.
  Then $f|  (f^{-1} (\Omega_j) \cap \hat{\Gamma})$
  is a branched analytic cover of  $\Omega_j$
   of some finite order, $j= 1,2$. Therefore,
  after possibly shrinking $\tau'$,
    we can choose
  a neighborhood $\cal U$ of $\tau$ in $\Bbb C^n$
  such that
  $f|  (f^{-1} (\Omega_j) \cap (\hat{\Gamma}\cap \cal U))$
  is a branched analytic cover of  $\Omega_j$
  of  order  $m_j \geq 0$  with $m_j$ {\em at most\/} equal to $1$.
 Hypothesis (a) implies that not both $m_j$ can be equal to $0$.
 If $m_1=1$ and $m_2=1$ then
$f^{-1} (\Omega_j) \cap (\hat{\Gamma}\cap \cal U)$ is
a graph of an analytic map $F_j$
on $\Omega_j$ for $j=1 $ and $j=2$.
The graphs $F_1$ and $F_2$  have identical boundary  values
on $\tau'$
equal  to $(f|\tau)^{-1}$ and  therefore continue
analytically  across $\tau'$
to give a single analytic  map $F$ on $\omega$. This
implies that $\hat{\Gamma} \cap \cal U$ is an analytic variety
and this means that $\tau$ is disjoint from the Shilov
boundary of ${\bf P}(\hat{\Gamma})$. This contradicts the
hypothesis  (b). Thus we are left only  with the case that
exactly one of the $m_j=1$ and the other multiplicity is $0$.
 Then the map $F_j$
extends smoothly to $\tau'$ and
 parametrizes $(\hat{\Gamma},\Gamma)$
near points of $\tau$ as a $2$-manifold with boundary. Therefore
$\tau$ is disjoint from $E$. This is a
contradiction and the  lemma follows. \hfill\qed

\nonumproclaim{Lemma 1.5} Let $\gamma$ be a finite union of
smooth curves in $\Bbb
C^n$  and let $x \in \hat{\gamma} \setminus \gamma${\rm .} There exists
a polynomial $P$ in $\Bbb C^n$ such that $P(x)=0$ and
$P \not = 0$ on $\hat{\gamma} \setminus \{x\}${\rm .}
\endproclaim

{\it Proof}.
We claim that there exists a complex linear map
$\phi=(\phi_1,\phi_2)\break :  \Bbb C^n \rightarrow \Bbb C^2$ such that
$(\phi|\hat{\gamma})^{-1}(\phi(x))= \{x\}$.
Set $V= \hat{\gamma} \setminus \gamma$; by [St],
$\hat{\gamma} \setminus \gamma$ is a $1$-dimensional subvariety of
$\Bbb C^n
\setminus \gamma$.
First choose a linear function
$\phi_1$ so that $\phi_1(x) \not\in \phi_1(\gamma)$.  Set
$q_1=\phi_1(x)$.
Then  $\phi_1^{-1}(q_1) \cap V$, being a $0$-variety bounded away
from
$\gamma$,  is a finite set
$\{y_1=x,y_2, \cdots,y_m\} \subseteq \Bbb C^n$. Choose a linear
function
$\phi_2$ such that $\phi_2$ separates the $m$ points
$\{y_1,y_2, \cdots,y_m\} \subseteq \Bbb C^n$.
Then  $\phi=(\phi_1,\phi_2) :  \Bbb C^n \rightarrow \Bbb C^2$
satisfies   $(\phi|\hat{\gamma})^{-1}(\phi(x))= \{x\}$.

Set $\gamma_0= \phi(\gamma) \subseteq \Bbb C^2$,
 $V_0=\phi(V) \subseteq \Bbb C^2$ and  $q= \phi(x) $.
 By the maximum principle, $V_0 \subseteq \widehat{\gamma_0}$;
also
$q \in \widehat{\gamma_0} \setminus \gamma_0$ and
$(\phi|\hat{\gamma})^{-1}( q)= \{x\}$.

Let $\ell$  be an affine complex line in $\Bbb C^2$ such that
$q $ is an isolated point  in $\ell \cap \widehat{\gamma_0}$. (Recall
that
$\widehat{\gamma_0} \setminus \gamma_0$ is a $1$-dimensional
subvariety
of $\Bbb C^2 \setminus \widehat{\gamma_0}$.)  Since
$\widehat{\gamma_0}$ is
polynomially convex, we can  find Runge domains in $\Bbb C^2$ that
decrease down to $\widehat{\gamma_0}$. In particular,  there exists
a Runge
domain $\Omega$ containing $\widehat{\gamma_0}$ such that,
$L$, the {\it connected component\/}  of
$\Omega \cap \ell $  that contains $q$,
satisfies $L \cap \widehat{\gamma_0} =\{q\}$ and
    $L$ is a hypersurface  in $\Omega$.
    By Serre [S] and Andreotti-Narasimhan[AN],
   since $\Omega$ is Runge in $\Bbb C^2$,
   $\check H^2(\Omega,\Bbb Z)=0$.  Hence by the Cousin II problem,
   there exists a function $F_0$ holomorphic on $\Omega$ such that
   $L=\{z \in \Omega: F_0(z)=0\}$. In particular, $F_0\neq 0$ on
   $\widehat{\gamma_0} \setminus \{q\}$
   and  $dF_0(q) \not = 0$.
   Set $F= F_0 \circ \phi$. Then, since $V_0 \subseteq
\widehat{\gamma_0}$,
   $F$ is a holomorphic function on a neighborhood of
   $\hat{\gamma}$  and the  only
   zero of $F$ on $\hat{\gamma}$ occurs  at $x$. Approximating
   $F$ uniformly   on a neighborhood of $\hat{\gamma}$ by
polynomials
   then gives
   the desired   $P$. \hfill\qed

\medbreak {\it {\rm C.} The Arens-Royden theorem}. Let $K$ be a compact space. We
denote the
algebra of continuous complex-valued functions on $K$ by ${\bf
C}(K)$ and
denote the invertible elements  (i.e.\ nonvanishing functions)
in ${\bf C}(K)$ by ${\bf C}^{-1}(K)$. Then ${\bf C}^{-1}(K)$ is an
abelian group
under multiplication and
contains the subgroup
$\mbox{exp}({\bf C}(K)) =\{e^f :  f \in {\bf C}(K)\}$. By a theorem of
Bruschlinsky
the quotient group  ${\bf C}^{-1}(K)/\mbox{exp}({\bf C}(K))$ is
naturally isomorphic to ${\check H}^1(K,\Bbb Z)$, the first { \u C}ech
cohomology group
with integer coefficients.

Let $A$ be a uniform algebra on $K$.
Denote the invertible elements
in $A$ by $A^{-1}$.
(If $K$ is the maximal ideal space of $A$ then
$A^{-1}$ is just the set of  $f \in A$
such that $f \not = 0$ on $K$.)
The multiplicative group  $A^{-1}$ contains the subgroup
$\mbox{exp}(A) =\{e^f :  f \in A\}$.
The Arens-Royden theorem states that
if $K$ is the maximal ideal space
of $A$, then
the quotient group  $A^{-1} /\mbox{exp}(A)$ is
naturally isomorphic to ${\check H}^1(K,\Bbb Z)$. Moreover, this isomorphism  factors
through
the previous one in the sense that the natural map
$A^{-1} /\mbox{exp}(A) \rightarrow  {\bf C}^{-1}(K)/\mbox{exp}({\bf
C}(K))$
induced by the inclusion  $A \rightarrow {\bf C}(K)$ is an
isomorphism.

For $K$ a compact subset
of $\Bbb C^n$,  if  $K$ is polynomially convex, then  $K$ is the
maximal
ideal space of ${\bf P}(K)$  and we have the natural isomorphism
$ j: {\bf P}^{-1}(K)/\mbox{exp}({\bf P}(K))\rightarrow  {\bf C}^{-
1}(K)/\mbox{exp}({\bf C}(K))$
provided by the Arens-Royden theorem.
(In this setting, an easy proof of the Arens-Royden theorem can be
obtained by  approximating $K$ by Runge domains
$\Omega$, applying the
fact  ([GR, Th.~7, p.~250]) that $\check H^1(\Omega, \Bbb Z)
  \simeq  {\check H}^0(\Omega, {\cal O}^*)/
\exp ({\check H}^0(\Omega, {\cal O}))$, and taking the inductive
limit over $\Omega$.) The isomorphism $j$
reduces the problem of
finding a polynomial on $K$ with certain periods to
producing  a nonvanishing continuous function with those periods.

\nonumproclaim{Lemma 1.6} Let $K$ be a  polynomially convex
compact subset of
$\Bbb C^n$ and let $\sigma$
be a $1$\/{\rm -}\/cycle contained in $K$. Let $f \in   {\bf C}^{-1}(K)${\rm .}
Then there exists
a polynomial $P$ in $\Bbb C^n$ such that
\smallbreak
\centerline{${\displaystyle\Delta_\sigma (\hbox{\rm arg }P)=\Delta_\sigma (\hbox{\rm  arg
}f).}$}
\endproclaim

 \demo{{R}emark}   The notation  $\Delta_\sigma (\hbox{\rm arg }f)$
denotes the variation of the argument of $f$ along the oriented $1$-cycle
$\sigma$.
If $f$ and $\sigma$ are smooth,
$i \Delta_\sigma (\hbox{\rm arg } f) =\int_\sigma df/f$. Alternatively,
$\Delta_\sigma (\hbox{\rm arg }f) /2\pi $ is the degree of
$f/|f|$ as a map  $\sigma \rightarrow
S^1$.
\enddemo

\demo{Proof}  Let $[f]$ be the class of $f$ in
${\bf C}^{-1}(K)/\mbox{exp}({\bf C}(K))$. Since the natural map
$ {\bf P}^{-1}(K)/\mbox{exp}({\bf P}(K))\rightarrow  {\bf C}^{-
1}(K)/\mbox{exp}({\bf C}(K))$
is surjective, there exists $F\in {\bf P}^{-1}(K)$ such that
$[F]=[f]$, where $[F] \in  {\bf C}^{-1}(K)/\mbox{exp}({\bf C}(K))$.
Hence there exists $u \in {\bf C}(K)$ such that  $F=f  e^u$.
Since $\Delta_\sigma (\hbox{\rm arg }e^u) =0$, we get
$\Delta_\sigma (\hbox{\rm arg }F)=\Delta_\sigma (\hbox{\rm arg }f)+
\Delta_\sigma (\hbox{\rm arg }e^u)= \Delta_\sigma (\hbox{\rm arg }f)$.
Finally we can approximate $F$ uniformly on $K$ by a polynomial
$P$ so that
$|\Delta_\sigma (\hbox{\rm arg }P)-\Delta_\sigma (\hbox{\rm arg }F)
|\break <2\pi$.
Therefore
$\Delta_\sigma (\hbox{\rm arg }P)=\Delta_\sigma (\hbox{\rm arg }F)   $,
since  $\sigma$, being  a cycle (``closed''),  each $\Delta_\sigma $ term
is an integral multiple of $2\pi$.
Hence $\Delta_\sigma (\hbox{\rm arg }P)=\Delta_\sigma (\hbox{\rm arg }F)\break=
\Delta_\sigma (\hbox{\rm arg }f)$.\enddemo

To apply   Lemma 1.6  we shall need
the following explicit version of the Bruschlinsky theorem. 

\nonumproclaim{Lemma  1.7} Let $S_0$ be a compact bordered
Riemann surface{\rm ,} not necessarily connected{\rm ,}
and let $F$ be a finite subset of $S_0${\rm .}
Let $S$  be obtained from $S_0$ by identifying  points in
the classes of some partition
of $F${\rm .}  Let $\gamma$ be a    disjoint union of
Jordan curves in $S$ such that
$[\gamma] \not = 0$ in $H_1(S,\Bbb Z)${\rm .} Then there exists
$f \in  {\bf C}^{-1}(S)$ such that
$\Delta_\gamma( \arg f) < 0${\rm .}
\endproclaim

\demo{Proof}
 We claim that $H_1(S,\Bbb Z)$ is torsion-free.
 Let $p: S_0 \rightarrow S$ be the identification map.
 We verify the claim  in three
 steps. (a) $H_1(S_0,F;\Bbb Z)$ is torsion-free.
 This follows from the exact sequence
 $$0 \rightarrow H_1(S_0,\Bbb Z)
  \rightarrow H_1(S_0,F,\Bbb Z) \rightarrow H_0(F,\Bbb Z)  $$
 and the fact that $H_1(S_0,\Bbb Z) $ is
 torsion-free, as is $H_0(F,\Bbb Z)$.
 (b) The induced map
  $p_*: H_1(S_0,F;\Bbb Z)  \rightarrow
  H_1(S,p(F );\Bbb Z)$ is an isomorphism,
 as is easily checked using   Mayer-Vietoris
 sequences to localize at $p(F)$
 and to separate the branches of $S$.  (c)
 From (a) and (b) we conclude that
 $H_1(S,p(F ); \Bbb Z)$ has no torsion and
 hence our claim follows from
 the exact sequence
 $$0 \rightarrow H_1(S,\Bbb Z) \rightarrow  H_1(S,p(F);\Bbb Z).$$

  We  write $\  [S,S^1]$ for  the set of homotopy equivalence
 classes of continuous functions
 $f: S \rightarrow S^1 \subseteq \Bbb C$. In this case homotopy
 equivalence is the same as equivalence mod $e^{{\bf C}(S)}$
 in ${\bf C}^{-1}(S)$.
 As $S$ is a CW complex (even a finite simplicial complex)
 we can apply a classification theorem (see Spanier, [Sp,  Th.~8.1.8,
 p.~427]) to  conclude that there is a natural isomorphism
 $\psi: [S,S^1] \rightarrow H^1(S, \Bbb Z)$ given by
 $\psi([f])([\beta]) = \Delta_\beta( \arg f)  $,
 for all  continuous functions
 $f: S \rightarrow S^1 \subseteq \Bbb C$ and all $1$-cycles $\beta$ in
$S$.
 Hence, for all
 $T \in  H^1(S, \Bbb Z) =\mbox{Hom}(H_1(S, \Bbb Z),\Bbb Z)$, there
exists
 $f \in {\bf C}^{-1}(S)$ such that for all $1$-cycles $\beta$ in $S$,
 $T([\beta])= \Delta_\beta( \arg f)  $.

 Finally since  $[\gamma] \not = 0$ in $H_1(S,\Bbb Z)$
  and since $H_1(S,\Bbb Z)$ is torsion-free, there exists
  $T \in H^1(S, \Bbb Z) =\hbox{\rm Hom}(H_1(S, \Bbb Z),\Bbb Z)$
  such that $T([ \gamma]) \not = 0$.  By the previous paragraph,
 there exists
 $f \in {\bf C}^{-1}(S)$ such that
 $  \Delta_ \gamma( \arg f)   = T([\gamma]) \not = 0$.
 If   $  \Delta_ \gamma( \arg f)    < 0$ we are done;  otherwise we
 replace $f$ by  $1/f$.  \enddemo

\section{Proof of Theorem 3}

   \nonumproclaim{Lemma 2.1} Let $\gamma$ be a smooth
   compact oriented $1$\/{\rm -}\/chain in
   $\Bbb C^n$ satisfying the moment condition{\rm .}
   Let
   $$T = \sum n_j [V_j]$$
   be the unique holomorphic $1$\/{\rm -}\/chain
   in $\Bbb C^n \setminus \gamma$
   such that $bT=[\gamma]$ whose
   existence is given by the Harvey\/{\rm -}\/Lawson theorem{\rm .}
   If $\gamma$ satisfies the linking condition{\rm ,} then
   $T$ is positive{\rm ;} i{\rm .}e{\rm .,}
   $n_j >0$
   for all $j${\rm .}
\endproclaim

   \demo{Proof}   Fix an index $k$ and  a point $x \in V_k$ such
   that $x \not\in \overline{V_j}$ for $j \not =k$. By the maximum
   principle,
   $\hbox{\rm supp }T  \subseteq\hat{\gamma}$.
   By Lemma 1.5 there exists a polynomial $P$
    in $\Bbb C^n$ such that $P(x)=0$ and
$P \not = 0$ on $\hat{\gamma} \setminus \{x\}$.
Hence  for $A= {\bf Z}(P)$,
$$0 \leq {\rm link}(\gamma,A)=\sum n_j  \cdot\#(  V_j,A).$$
   For $j \not = k$, $A \cap V_j = \emptyset$ and so
   $\#(V_j,A)=0$.  We have therefore $0 \leq  n_k \cdot\#(V_k,A) $.
   As $P(x)=0$,  $\#(V_k,A) >0$ and we get that
   $0 \leq n_k $. We conclude that $0 < n_k$.\enddemo

We first prove Theorem 3 in two special cases. 
\noindent

\nonumproclaim{{\it Case}  {\rm (i)}} $\gamma$ is a simple closed oriented smooth curve{\rm .}\endproclaim

 \demo{{P}roof} If $\gamma$ is polynomially convex, then
${\bf P}(\gamma)={\bf C}(\gamma)$. Hence, first
choosing an $f \in {\bf C}(\gamma)$ such that
$f$ maps $\gamma$ to the unit circle with
$\frac{1}{2\pi}\Delta_\gamma( \arg f) =-1$,
 we get a polynomial $P$ such that
$\frac{1}{2\pi}\Delta_\gamma( \arg P) =-1$.
Therefore, by Lemma 1.2,   ${\rm link}(\gamma,A)= -1$
where $A={\bf Z}(P)$, and this
contradicts the linking condition.
We conclude that $\gamma$ is not polynomially convex.
It follows that
$V= \hat{\gamma} \setminus \gamma$ is a 1-variety
of finite area in $\Bbb C^n
\setminus \gamma$ and $b[V]=[\gamma]$, as currents;
cf. Lemma 2.4  below.
Let $\psi$ be a holomorphic $(1,0)$-form in $\Bbb C^n$.
Then, since $[V]$ is a $(1,1)$-current and $ d \psi$ is a $(2,0)$-form, we
get
$$\int_\gamma \psi =  [V]( d\psi)=0.$$
 This says that $\gamma$ satisfies the moment condition,  proving
 the theorem in case (i). \enddemo

\nonumproclaim{{\it Case} {\rm (ii)}} $\gamma$ is a real analytic $1$\/{\rm -}\/cycle{\rm .}
\endproclaim

 \demo{{P}roof}
We claim that  $[\gamma]  = 0$ in $H_1(\hat{\gamma},\Bbb Z)$.
Suppose, by way of contradiction,  
that $[\gamma] \neq 0$ in $H_1(\hat{\gamma},\Bbb Z)$.
Since $\hat{\gamma}$ is a compact bordered Riemann surface with
 a finite number of
points identified, we can apply Lemma 1.7 to obtain
an $f \in {\bf C}^{-1}(\hat{\gamma})$ such that
$\frac{1}{2\pi}\Delta_{\gamma} \hbox{\rm arg } f  < 0$.
By the Arens-Royden theorem
in the form of Lemma 1.6,  there is  a polynomial $P$
such that
 $\frac{1}{2\pi}\Delta_{\gamma} \hbox{\rm arg } P
 = \frac{1}{2\pi}\Delta_{\gamma} \hbox{\rm arg } f <0$.
 This contradicts the linking condition. We conclude that
 $\gamma  \sim 0$ in $\hat{\gamma}$.

 We claim that $\gamma $  satisfies the moment   condition.
 Then,  by Lemma 2.1, Theorem 3    follows in this case.
 Let $\psi$ be a holomorphic 1-form in $\Bbb C^n$.
 Then $d\psi$ is a $(2,0)$-form and so
 $d\psi=0$ on the one-dimensional analytic set $\hat{\gamma}$.
  But, since  $\gamma  \sim 0$ in $\hat{\gamma}$,
  $\gamma =b\Sigma$ where $\Sigma$ is a $2$-chain in
$\hat{\gamma}$.
  Hence by Stokes' theorem,
  $\int_\gamma \psi = \int_\Sigma d\psi = 0$. This is the moment
  condition.\enddemo

\nonumproclaim{{\it Case} {\rm (iii)}} The general case{\rm .}\endproclaim

 \demo{{P}roof} Arguing as in case (i) we see that $\gamma$
 cannot be polynomially convex. Thus we can suppose
  that  $\gamma $ is not polynomially convex in the general case. Choose
a minimal subfamily ${\cal F} \subseteq \{   \gamma_{j}  \}$ such
that the
polynomial hull of the sum  $\Gamma  =\sum \{\gamma_j :
\gamma_j \in
{\cal F}\}$
satisfies
 $\hat{\Gamma} \setminus \gamma = \hat{\gamma} \setminus
\gamma$;
 then ${\cal F} \not = \emptyset$
 because $\gamma $ is not polynomially convex.
Let  $\sigma  =\sum \{\gamma_j : \gamma_j  \not\in
{\cal F}\}$. We get a partition of $\gamma$ as
$\gamma= \Gamma +\sigma$.
We are abusing language somewhat, since we write $\gamma$,
$\Gamma$
and $\sigma$ as
oriented $1$-cycles and also,
when we take the    polynomially convex hull,  as the
corresponding underlying sets  in $\Bbb C^n$.
Let $V = \hat{\Gamma} \setminus \Gamma$; $V$ is a one
dimensional subvariety
of $\Bbb C^n  \setminus \Gamma$.
Let  ${\frak S}$ denote the Shilov boundary of  the uniform
algebra  ${\bf P}(\hat{\Gamma})$. \enddemo

\nonumproclaim{Lemma  2.2} {\rm (a)} $\Gamma  =  \frak S$  and 
 {\rm (b)}  $\Gamma \subseteq
\overline{\hat{\Gamma} \setminus \Gamma}${\rm .}\endproclaim

\demo{Proof} (a) Clearly   ${\frak S} \subseteq \Gamma$. We need only
show that
  $\Gamma  \subseteq  \frak S.$
Arguing by contradiction, we suppose otherwise. Then there exists an
open subarc $\tau$ of some $\gamma_k \in \cal F$ such that
$\tau \subset \widehat{\Gamma \setminus \tau}$.
Put ${\cal F}_1={\cal F} \setminus \gamma_k$,
 $\Gamma_1=\sum \{\gamma_j : \gamma_j   \in
{\cal F}_1\}$ and $\beta= \gamma_k \setminus \tau$.
Then
$ \Gamma \setminus \tau = \Gamma_1 \cup \beta   $.
By Lemma 1.3,  $\widehat{(\Gamma_1 \cup \beta})=
\widehat{\Gamma_1} \cup \beta$
and so $\widehat{\Gamma_1} \cup \beta$ is polynomially convex.
Hence $\Gamma \subseteq  \widehat{\Gamma_1} \cup \beta$.
Therefore $\hat{\Gamma} \setminus \gamma=
\widehat{\Gamma_1} \setminus \gamma$. Thus, as
${\cal F}_1$  is a proper subset of $\cal F$,
this contradicts the minimality of $\cal F$.
Part (a)  follows.

 (b)  We argue by contradiction and suppose that
there exists an
open subarc $\tau$ of some $\gamma_k \in \cal F$ such that
$\tau$ is disjoint from
$  \overline{\hat{\Gamma} \setminus \Gamma}$.  Then, by the
local maximum modulus principle,
$\hat{\Gamma} \setminus \Gamma \subseteq
\widehat{\Gamma \setminus \tau}$.
As in part  (a),
$ \Gamma \setminus \tau = \Gamma_1 \cup \beta   $
and $\hat{\Gamma} \setminus \gamma=
\widehat{\Gamma_1} \setminus \gamma$.
Again this contradicts the minimality of $\cal F$
and part (b)  follows. \enddemo

The next lemma is due  essentially  to Lawrence [L], who treats the
case of
 a  simple closed rectifiable curve $\Gamma$. We shall briefly
indicate how his proof adapts to our setting, in which $\Gamma$ is
smooth, but not connected.

\nonumproclaim{Lemma 2.3} The $1$\/{\rm -}\/variety
$V= \hat{\Gamma} \setminus \Gamma $ has  finite
area {\rm (}${\cal H}^2$  measure\/{\rm )}  and  the corresponding  positive
$(1,1)$\/{\rm -}\/current $[V]$ {\rm (}\/oriented by the natural orientation of $V${\rm )}
satisfies
$$  b[V] = \sum \{\varepsilon_j[\gamma_j] : \gamma_j \in
{\cal F}\},\leqno(2.1)$$
where each $\varepsilon_j=\pm 1${\rm .}\endproclaim

 \demo{{R}emark} We do not use the linking hypothesis in Lemma 2.3.
 Without that
hypothesis, it is, in general,  not true that $b[V] =  [\Gamma] $
for  $V = \hat{\Gamma} \setminus \Gamma$ with
$\Gamma$ a $1$-cycle in $\Bbb C^n$.  For example, take
$\Gamma$ to be the unit circle in $\Bbb C$ with the
clock-wise orientation; then $V$ is the open unit disk and
$b[V] =  -[\Gamma] $. It is the
addition of the linking hypothesis for $\gamma$ that will
yield the correct orientation in Lemma 2.4.\enddemo

\demo{Proof} Lawrence's  argument that $\hat{\Gamma} \setminus
\Gamma $
has finite area is  valid when $\Gamma$ is a finite
 union of simple closed smooth curves.
 Hence the $(1,1)$-current   $[V]$ exists
 with $\hbox{\rm supp}([V] \subseteq \Gamma$.
 Lawrence's arguments, together
 with  Lemma 2.2,  imply that
 $b[V]= {\cal H}^1\LL \Gamma \wedge \eta$ where $\eta $ is a
 Borel measurable {\em unit\/}
 tangent vectorfield
 to $\Gamma$; in particular
 $b[V]$ has multiplicity $1$ at almost every point of $\Gamma$.
 Finally Lawrence's argument shows that
 $(b[V]) \LL  \gamma_j = \pm [\gamma_j]$ for each $\gamma_j \in
\cal F$,
 since $[\gamma_j]$ is an indecomposable
   integral current. This gives (2.1).\enddemo

\nonumproclaim{Lemma  2.4}  With $V$ as in Lemma {\rm 2.3,}
   $$b[V] =  [\Gamma] .$$
\endproclaim

\demo{Proof}  We need to show that   $\varepsilon_j=1$ for all $j$.
 Fix an index $k$ with $\gamma_k \in {\cal F}$.
 Since $\gamma_k \subseteq \frak S$ by Lemma 2.2 (a),
  we can choose a polynomial
 $F$ so that  $F(x)=1$ for some $x \in \gamma_k$ and
 $ |F| < 1/2$ on  the set $\gamma \setminus \gamma_k$. Choose,
 by Lemma 2.2 (b),  a point
 $\lambda \in  F(\hat{\gamma}) \setminus F(\gamma)$
 with   $|\lambda| > 1/2$ and
 set $A= {\bf Z}(F- \lambda)$, a complex hypersurface in $\Bbb C^n$.
Then $A$  is disjoint from $\gamma$ and,
by the linking hypothesis on $\gamma$,
${\rm link}(\gamma,A) \geq 0$.

On  all   $\gamma_j$, $j \not = k$,
 $ |F| < 1/2 < |\lambda|$; hence $F-\lambda$ has a logarithm on
$\gamma_j$
 and so
 $$  \frac{1}{2\pi i}\int_{\gamma_j}
 \frac{d(F-\lambda)}{F-\lambda}=0.\leqno(2.2)$$
Hence
$$\frac{1}{2\pi i}\int _\gamma \frac{d(F-\lambda)}{F-\lambda}=
\sum_j \frac{1}{2\pi i}\int_{\gamma_j}  \frac{d(F-\lambda)}{F-
\lambda}=
\frac{1}{2\pi i}\int_{\gamma_k}  \frac{d(F-\lambda)}{F-\lambda}.$$
 Therefore  we have
 $$  0\leq  {\rm link}(\gamma,A) =
 \frac{1}{2\pi i}\int _\gamma \frac{d(F-\lambda)}{F-\lambda}=
 \frac{1}{2\pi i}\int_{\gamma_k}  \frac{d(F-\lambda)}{F-\lambda}.
\leqno(2.3)$$
 From this we will deduce that $\varepsilon_k >0$. We can assume that
 $\lambda$ is a regular value of $F|V$ and that
 $V$ is $s$-sheeted over the component
 $\Omega$ of $\Bbb C \setminus F(\gamma)$
 containing $\lambda$, $s \geq 1$. Hence
 there exists a small
 closed disk $\Delta \subseteq \Omega$
  centered at $\lambda$  such that $F^{-1}(\Delta) \cap V$ is the
disjoint
  union of $s$ components $\Delta_1,\Delta_2, \cdots,\Delta_s$ each
of
  which is mapped biholomorphically to $\Delta$.
  By (2.1),
  $$b[V \setminus \cup_{i=1}^s \Delta_i]=b[V]-\sum_{i=1}^s
b[\Delta_i]=
  \sum \{\varepsilon_j[\gamma_j] : \gamma_j \in
{\cal F}\}-\sum_{i=1}^s b[\Delta_i].$$
Hence we get, as
 $\omega=\frac{1}{2\pi i}\frac{d(F-\lambda)}{F-\lambda}$ is a
closed
$1$-form on $V \setminus \cup_{i=1}^s \Delta_i$,
\begin{eqnarray*}
0=[V \setminus \cup_{i=1}^s \Delta_i](d\omega)&=&
b[V \setminus \cup_{i=1}^s \Delta_i](\omega)\\
&=&
  \sum \{\varepsilon_j[\gamma_j] (\omega): \gamma_j \in
{\cal F}\}-\sum_{i=1}^s b[\Delta_i](\omega). 
\end{eqnarray*}
By the Cauchy integral formula,
$$b[\Delta_i](\omega)=\frac{1}{2\pi i} \int_{b\Delta}\frac{dz}{z-
\lambda}=1$$
for each index $i$.
Thus, applying   (2.2),  we get
$$0=
\varepsilon_k\frac{1}{2\pi i}\int_{\gamma_k}  \frac{d(F-\lambda)}{F-
\lambda}
-s.\leqno(2.4)$$
 Since $s\not = 0$, (2.4) implies that
 $$\frac{1}{2\pi i}\int_{\gamma_k}  \frac{d(F-\lambda)}{F-
\lambda}\not
 =0.$$
 Hence, by (2.3),
 $$\frac{1}{2\pi i}\int_{\gamma_k}  \frac{d(F-\lambda)}{F-\lambda}
>0.$$
 Now (2.4) implies that $\varepsilon_k >0$. Therefore $\varepsilon_k=1$
and this
gives the lemma. \enddemo

   Now consider the above  partition
 $\gamma$ as $\gamma = \Gamma +\sigma$, where
 $\hat{\Gamma} \setminus  \gamma =\hat{\gamma} \setminus
\gamma$.
 We set  $V=  \hat{\Gamma} \setminus  \Gamma $.
 Consider the two cases:
\medbreak 1.  $\sigma  \not\subseteq  V$, 
 or
 2.  $\sigma   \subseteq  V$. 

\demo{{C}ase {\rm 1}}    Fix  $x \in  \sigma$ with $x \not\in \hat{\Gamma} $.  Then
 $x \in \gamma_k$ for some $\gamma_k$ which is not one of the
curves
 which comprise $\Gamma$. We construct a smooth
 complex-valued function $f$ on  $\hat{\gamma}$
 as follows: first take
 $f \equiv 1$ on all of  $\hat{\gamma}$ except for a small
 subarc $v$ of $\gamma_k$ such that
 $x \in v$ and $\bar{v} \cap  \hat{\Gamma} =  \emptyset$.
  We can  then  extend  $f$ so that the
 image of $f$ on $v$ winds once negatively about the unit circle.
 Then $f$ is   nonvanishing on $\hat{\gamma}$. As $f \equiv 1$ on
 $\hat{\Gamma}$, we have $f \in {\bf P}(\hat{\Gamma})$.
 By the hypothesis for case 1,  $\hat{\gamma} $ is the union of
 $\hat{\Gamma}$ and some of the $\sigma$
  curves which are not contained in
 $\hat{\Gamma}$.
 Hence (see [St])
 $f \in {\bf P}(\hat{\gamma})$.  By our
 construction
 $$\frac{1}{2\pi i}\int_{\gamma_k}  \frac{df}{f}=
 \frac{1}{2\pi i}\int_{v}  \frac{df}{f}=-1$$
 since $f \equiv 1$ on $\gamma_k \setminus v$,
 and
  $$\frac{1}{2\pi i}\int_{\gamma_j} \frac{df}{f}= 0$$
for $j \not = k$, since $f \equiv 1$ on these $\gamma_j$.
Hence
$$\frac{1}{2\pi i}\int_{\gamma} \frac{df}{f}= -1.$$
Approximating $f \in {\bf P}(\hat{\gamma})$ we get a polynomial
$P$
such that $P \not = 0$ on $\hat{\gamma}$ and
$$\frac{1}{2\pi i}\int_{\gamma} \frac{dP}{P}= -1.$$

By Lemma 1.6, this gives
 ${\rm link}(\gamma,A)=-1$ where  $A=\{z \in \Bbb C^3:P(z)=0\}$.
This contradicts the linking condition and
we conclude that Case 1 cannot arise.\enddemo

\demo{{C}ase {\rm 2}}  Let $E$ be the set of ``bad'' points of $\Gamma$ given
in Lemma 1.4.
Locally at
each point of $\Gamma \setminus E$,
$\hat{\Gamma}\setminus \Gamma$ is a $2$-manifold with boundary.
 Choose $\psi$ a real-valued
${\cal C}^\infty$ function
on $\Bbb C^n$ such that $\psi \geq 0$ on $\Bbb C^n$ and
$\psi=0$ on $E$. Choose, by Sard's theorem,  $\varepsilon >0$
(``admissible'') so that
\begin{itemize}
\item[(a)] $\varepsilon$ is a regular value of $\psi| \Gamma$, 
\item[(b)] $\varepsilon$ is a regular value of  $\psi| V_{\rm reg}$,
\item[(c)] $\psi \neq \varepsilon$ on $V_{\rm sing }$.
\end{itemize}

Set $D_\varepsilon=\{z \in \hat{\Gamma}: \psi \geq \varepsilon\}$,
$Q_\varepsilon=\{z \in \hat{\Gamma}: \psi  \leq \varepsilon\}$ and
$\alpha_\varepsilon =\{z \in \hat{\Gamma}: \psi = \varepsilon\}$.
Then $\alpha_\varepsilon= \tau_\varepsilon + \rho_\varepsilon$ where
$\tau_\varepsilon$ is a finite  set of closed curves in $V_{\rm reg}$ and
$\rho_\varepsilon$ is a finite union of arcs joining two points of
$\Gamma$ and, except for its endpoints, lying in $V_{\rm reg}$.
Except for the  finite set $V_{\rm sing} \cap D_\varepsilon$, $D_\varepsilon$
is a topological manifold
 with boundary $bD_\varepsilon$, where $bD_\varepsilon$ is
piecewise smooth consisting of the
oriented curves $\tau_\varepsilon$ and other oriented
curves, whose sum   we denote by $\kappa_\varepsilon$;
thus $\kappa_\varepsilon$ is  a sum of some
   subarcs of $\Gamma$-curves and the arcs
of $\rho_\varepsilon$.  Thus $bD_\varepsilon=
\tau_\varepsilon+\kappa_\varepsilon$.

We consider two  subcases:
\medbreak
 {\it Case} 2a:  There exists an (admissible) $\varepsilon > 0$,
such that  $[\sigma] \neq  0$
in $H_1(D_\varepsilon,\tau_\varepsilon;\Bbb Z)$ 
or
\smallbreak

{\it Case} 2b:  For all (admissible) $\varepsilon > 0$, $[\sigma] = 0$
in $H_1(D_\varepsilon,\tau_\varepsilon;\Bbb Z)$.
\medbreak

Before considering case 2a we need two lemmas.
For  a nonvanishing
   continuous complex-valued function $h$  defined on
   an oriented  $1$-cycle $C$,  we denote the index of $h$ on $C$
   by $\hbox{\rm Ind}(h, C)$.  This equals both
   $\frac{1}{2\pi} \Delta_{C} (\hbox{\rm arg } h)$  and
   the winding number of the curve $h(C)$ about the origin.\enddemo

\nonumproclaim{Lemma  2.5}  Let $A$ be the planar annulus
   $\{z \in \Bbb C:  a \leq  |z| \leq b \}${\rm ,} $a <b${\rm ,} and let
   $\Gamma_a= \{z \in \Bbb C:    |z| =a \} $ and
   $\Gamma_b= \{z \in \Bbb C:    |z| =b \} ${\rm ,} both positively
oriented{\rm .}
\begin{itemize}
   \ritem{(a)} Let $h$ be a nonvanishing
   continuous complex\/{\rm -}\/valued function defined on
   $bA= \Gamma_b - \Gamma_a$ such that
   $\hbox{\rm Ind}(h,\Gamma_a) =   \hbox{\rm Ind}(h,\Gamma_b)${\rm .}
   Then  $h$ extends to be a  nonvanishing continuous
   complex\/{\rm -}\/valued function
   $H$ defined on
   $A${\rm .}
\ritem{(b)}  Let $g$ be a nonvanishing continuous
   complex\/{\rm -}\/valued function defined on
   $ \Gamma_b \cup   S$  where $S$ is a proper closed subset of
   $\Gamma_a${\rm .} Then  $g$ extends
   to be a  nonvanishing continuous
   complex\/{\rm -}\/valued function defined on
   $A${\rm .}
\end{itemize}

\endproclaim
 
 \demo{Proof} (a) This is a special case of a more more general extension
   theorem of Hopf. We  give a short proof for our special case. Let
   $k= \hbox{\rm Ind}(h,\Gamma_a) \in \Bbb Z$. Set $q= h z^{-k}$ on
   $bA$. Then $q$ satisfies
   $\hbox{\rm Ind}(q,\Gamma_a) =   \hbox{\rm Ind}(q,\Gamma_b)= 0$.
   Hence $q $ has a complex-valued logarithm $u$ on $bA$, i.e.
   $q=e^u$ on $bA$. By Tietze's extension theorem, we can extend $u$
to be
   a continuous complex-valued function $\Phi$ on $A$.
   Now take $H = z^{k} e^\Phi$ on $A$. 

  (b) By part (a) it suffices to extend $g$ to a nonvanishing
continuous
  complex-valued function on  $\Gamma_a$ such that
  $\hbox{\rm Ind}(g,\Gamma_a) =   \hbox{\rm Ind}(g,\Gamma_b)$.
  Let $\tau$ be an open subarc of   $\Gamma_a$
  whose closure is disjoint from $S$.  Then  $\Gamma_a  \setminus
\tau$
  is a  closed interval  containing $S$.  Hence we can extend $g$ from
  $S$ to be a continuous complex-valued nonvanishing function on
  $\Gamma_a \setminus \tau$ with $g=1$ on the two  endpoints of
  $\Gamma_a  \setminus \tau$. Since
   $g=1$ on the two  endpoints of
  $\Gamma_a \setminus \tau$,
  $ \frac{1}{2\pi}\Delta_{   \Gamma_a  \setminus \tau  }  (\hbox{\rm arg
} g)
  \in \Bbb Z.$
  Hence $j=\hbox{\rm Ind}(g,\Gamma_b) -
   \frac{1}{2\pi}\Delta_{   \Gamma_a  \setminus \tau  }  (\hbox{\rm arg }
g)
   \in \Bbb Z$.
   Now  we can extend $g$
  over $\tau$  such that, on $\tau$, $g$ is
   a map  covering the unit circle   $j$ times; that is,
  $g$ on $\tau$ has  complex values  of modulus one   and satisfies
   $\frac{1}{2\pi} \Delta_{\tau} (\hbox{\rm arg } g)=  j$.
    Thus
    $\hbox{\rm Ind}(g,\Gamma_a) =
    \frac{1}{2\pi} \Delta_{\tau} (\hbox{\rm arg } g) +
    \frac{1}{2\pi} \Delta_{   \Gamma_a  \setminus \tau  }  (\hbox{\rm arg
} g)
    = \hbox{\rm Ind}(g,\Gamma_b)$,
   as desired. \enddemo

  \nonumproclaim{Lemma  2.6} Let $A$ be the planar annulus
   $\{z \in \Bbb C:  a \leq  |z| \leq b \}${\rm ,} $a<b${\rm ,}  and let
   $\Gamma_a= \{z \in \Bbb C:    |z| =a \} $ and
   $\Gamma_b= \{z \in \Bbb C:    |z| =b \} ${\rm ,} both positively
oriented{\rm .}
\begin{itemize}
   \ritem{(a)} Let $h$ be a nonvanishing continuous
   complex\/{\rm -}\/valued function defined on
    $|z| \leq b${\rm .} Then there exists
    a  continuous  complex\/{\rm -}\/valued function  $f$ defined on
    $|z| \leq b$ such that
\begin{itemize}
   \ritem{i)} $f \neq 0$ on $|z| \leq b${\rm ,}
    \ritem{ii)} $f=1$ on $|z| \leq a${\rm ,} and 
    \ritem{iii)} $f=h$ on $\Gamma_b${\rm .}
\end{itemize}
    \ritem{(b)} Let $h$ be a nonvanishing continuous
    complex\/{\rm -}\/valued function defined on
    $ A${\rm .}   Let $S$ be a proper subset of $\Gamma_a${\rm .}
    Then there exists
    a  continuous  complex\/{\rm -}\/valued function  $f$ defined on
    $A$ such that
\begin{itemize}
    \ritem{i)} $f \neq 0$ on $A${\rm ,}
\ritem{ii)} $f=1$ on $S${\rm ,} and
\ritem{iii)} $f=h$ on $\Gamma_b${\rm .}\end{itemize}
\end{itemize}

\endproclaim
 
    \demo{Proof} (a) Set $f=1$ on $|z| \leq a$.  Since $\hbox{\rm ind}(h,
    \Gamma_b)=0$, we can use (a) of  Lemma 2.5  to extend $f$
    to $A$ so that $f = h$ on $\Gamma_b$.

    (b)  Define  $g$ on $\Gamma_b \cup S$ by
    $g \equiv 1$ on $S$ and $g= h$ on $\Gamma_b$. Then apply (b) of
     Lemma  2.5
    to extend $g$  as the desired function on $A$.
\enddemo

We first consider Case 2a: We assume
that  $[\sigma] \neq  0$
in $H_1(D_\varepsilon,\tau_\varepsilon;\Bbb Z)$
the first relative singular homology group.
Let $\tilde{D}_\varepsilon $ be obtained from $D_\varepsilon$ by attaching
a disk
at each of the components of $\tau_\varepsilon$.
Then $\tilde{D}_\varepsilon $ is a finite simplicial complex.
Let $\tilde{\tau}_\varepsilon$
be the union of these closed disks in  $\tilde{D}_\varepsilon $. There
is a natural map
$$H_1(D_\varepsilon,\tau_\varepsilon;\Bbb Z) \rightarrow
 H_1(\tilde{D}_\varepsilon,\tilde{\tau}_\varepsilon;\Bbb Z) $$
 induced by the inclusion $(D_\varepsilon,\tau_\varepsilon)\subseteq
 (\tilde{D}_\varepsilon,\tilde{\tau}_\varepsilon)$ and, by excision, this map
is
 an isomorphism. Hence
 $[\sigma] \neq  0$
in $H_1(\tilde{D}_\varepsilon,\tilde{\tau}_\varepsilon;\Bbb Z)$.
From the exact homology  sequence
$$ H_1(\tilde{D}_\varepsilon;\Bbb Z)
\rightarrow
H_1(\tilde{D}_\varepsilon,\tilde{\tau}_\varepsilon;\Bbb Z) $$
we conclude that
 $[\sigma] \neq  0$
in $H_1(\tilde{D}_\varepsilon ;\Bbb Z)$.
Since $\tilde{D}_\varepsilon$ is a compact bordered Riemann surface
with
 a finite number of
points identified, we can apply Lemma 1.7 to obtain
an $f \in {\bf C}^{-1}( \tilde{D}_\varepsilon)$ such that
$\frac{1}{2\pi}\Delta_{\sigma} \hbox{\rm arg } f  < 0$.
By
  Lemma  2.6,
since  $\tilde{\tau}_\varepsilon$ is a disjoint union of closed disks
disjoint from $\sigma$,
we can
arrange that $f \equiv 1$ on $\tilde{\tau}_\varepsilon$.
Since $\sigma \subseteq D_\varepsilon \subseteq \tilde{D}_\varepsilon$,
by restricting $f$ to $D_\varepsilon$ we get
$f \in  {\bf C}^{-1}(D_\varepsilon)$ such that
$\Delta_\sigma (\hbox{\rm arg }f  )< 0$ and
$f \equiv 1$ on $\tau_\varepsilon= D_\varepsilon \cap
\tilde{\tau}_\varepsilon$.

Since $\kappa_\varepsilon \cap
\rho_\varepsilon$ is a {\it proper\/} subset
of  each closed curve in $\kappa_\varepsilon$,
we can apply
   Lemma 2.6
to arrange that $f \equiv 1$  on $\rho_\varepsilon$.  Thus
$f \equiv 1$ on $ \tau_\varepsilon +\rho_\varepsilon
=\alpha_\varepsilon =bQ_\varepsilon$ and thus we can extend $f$  to be a continuous function on
$\hat{\Gamma}$  by setting $f \equiv 1$ on $Q_\varepsilon$.
 To summarize:
we have
$f \in {\bf C}^{-1}(\hat{\Gamma})$ such that
$\Delta_\sigma (\hbox{\rm arg }f  )< 0$.  Applying the Arens-Royden
theorem we get a polynomial $P$ such that $P \neq 0$ on
$\hat{\Gamma}$
and $\Delta_\sigma (\hbox{\rm arg }P)< 0$.  This contradicts the linking
condition. 

We next consider Case 2b.
We assume that for all (admissible) $\varepsilon > 0$, $[\sigma] = 0$
in $H_1(D_\varepsilon,\tau_\varepsilon;\Bbb Z)$.
We construct  a compact subset $K$
of $\hat{\Gamma} \setminus E$ by adjoining to the set $\sigma$,
for each component of $\sigma$,  a path in
$\hat{\Gamma} \setminus E$ joining that component to some point
of
$\Gamma \setminus E$.
If necessary we can further enlarge $K$ so that each of the finitely
many components of $\hat{\Gamma} \setminus \Gamma$
contains an arc joining one endpoint  to a point of $\Gamma
\setminus
E$.
We choose a sequence
of admissible ${\varepsilon_k}$ decreasing to $0$ and we write
$D_k$ for $D_{\varepsilon_k}$,  $Q_k$ for $Q_{\varepsilon_k}$, etc.
We choose $\varepsilon_1$ small enough so that
$K \subseteq  D_k$ for all $k$.

Then, since $[\sigma] = 0$
in $H_1(D_k,\tau_k;\Bbb Z)$
for all $k$, there exists a $2$-chain
$\Sigma_k$ in $D_k$ and a $1$-cycle $\mu_k$ in $\tau_k$ such that
$b\Sigma_k=\sigma + \mu_k$. For $x \in D_k \setminus ( \tau_k
\cup
\sigma)$, the multiplicity of $\Sigma_k$ at $x$, defined as the
intersection
number  of homology classes $[\Sigma_k]$ and $[x]$, is constant on
each
component of
$D_k \setminus ( \tau_k \cup  \sigma)$.
We write  $\hbox{\rm mult}(\Sigma_k,x)$ for
the multiplicity of $\Sigma_k$ at $x \in D_k \setminus ( \tau_k \cup
\sigma)$.
Let $T_k$ be the (1,1)-current associated to $\Sigma_k$.
Then $T_k$ is just integration over the finite number of components
of
$D_k \setminus ( \tau_k \cup
\sigma)$ with integral weights given by
the multiplicity of $\Sigma_k$. We have
$bT_k= [\sigma]+ [\mu_k]$.

Let $N_k$ be the maximum of
$|\hbox{\rm mult}(\Sigma_k,x)|$
(as $x$ varies  over the finite number of connected
components of
$D_k \setminus ( \tau_k \cup \sigma)   $ ).

\nonumproclaim{Lemma 2.7} $N_k \leq  N_1$ for all $k${\rm .}
\endproclaim

\demo{Proof}  We can assume that
the $2$-cycle $\Sigma_k$ is contained in
the union of those components
of $D_k$ which meet $\sigma $, since homology is the sum of the
homology
of the (path) components. Fix  $x \in D_k \setminus  ( \tau_1 \cup
\sigma) $
 such that $\hbox{\rm mult}(\Sigma_k,x) \neq 0$.
We shall show that
$|\hbox{\rm mult}(\Sigma_k,x)| \leq N_1$.

Let $\Omega$ be the
component of $D_k $ which contains $x$.
Then $\Omega$ meets $\sigma$ and
so contains   components of $\sigma$.
Let $\delta$ be an arc in $\Omega$ joining $x$ to a component
of $\sigma$ and otherwise disjoint from $\sigma$.
Then that component of  $\sigma$
and a path in
$\hat{\Gamma} \setminus E$ joining that component to some point
of
$p \in \Gamma \setminus E$
are  contained in a component $\Omega_1$ of $D_1 \setminus
\tau_1$  with
$\Omega_1 \subseteq \Omega$.

 Since $b(\Sigma_k -\Sigma_1) = \mu_k-\mu_1$ is disjoint from
 $\Omega_1$, the multiplicity of  $\Sigma_k -\Sigma_1$ in
 $\Omega_1$ is constant.  We claim
 that this multiplicity is $0$. In fact,  near $p \in \Gamma \setminus
E$,
  there is a
 deformation
 retraction of $\Omega_1$ to a smaller set which is disjoint from a
small
 neighborhood  $\omega$ of $p$. This retraction moves $\Sigma_k$
and
 $\Sigma_1$ off of $\omega$ but does not change
$b\Sigma_k$  or
 $b\Sigma_1$  and so does not change
 the multiplicity of  $\Sigma_k - \Sigma_1$
 on $\Omega_1$. As   the  modified 2-chains
 each have multiplicity $0$   on
 $\omega$ we conclude that $\Sigma_k   - \Sigma_1$
 has multiplicity  $0$ on $\Omega_1$. This gives the claim. We have
therefore
 $$\hbox{\rm mult}(\Sigma_k,y)=\hbox{\rm mult}(\Sigma_1,y)$$
 for each $y \in \Omega_1$.

  Since the multiplicity of $\Sigma_k$ is constant on
  components  of
   $\Omega \setminus \sigma$   we  have
 $$\hbox{\rm mult}(\Sigma_k,x)=\hbox{\rm mult}(\Sigma_k,y_0),$$
 for some $y_0 \in \Omega_1$--- we can just choose a point
 $y_0$  on
 $\delta  \cap \Omega_1$.
 We conclude that
 $$\hbox{\rm mult}(\Sigma_k,x)=\hbox{\rm mult}(\Sigma_1,y_0).$$
Thus $|\hbox{\rm mult}(\Sigma_k,x)|  \leq
|\hbox{\rm mult}(\Sigma_1,y_0)|   \leq N_1$.  \enddemo

 We denote the mass norm of a current $T$ by
${\bf M}(T)$. In particular,  ${\bf M}([Q_\varepsilon]) $
is then the ``area''  (${\cal H}^2$  measure) of $Q_\varepsilon$.

\nonumproclaim{Lemma 2.8} ${\bf M}([Q_\varepsilon ])\rightarrow 0$
as $\varepsilon \rightarrow 0${\rm .}\endproclaim

\demo{Proof}  $Q_\varepsilon \subseteq \hat{\Gamma}$ and ${\cal
H}^2(\hat{\Gamma} )
< \infty$. Since $Q_\varepsilon \rightarrow E$ as
$\varepsilon \rightarrow 0$ and ${\cal H}^2(E)=0$, we conclude that
${\cal H}^2( Q_\varepsilon)   \rightarrow 0$ as
$\varepsilon \rightarrow 0$.  \enddemo

\nonumproclaim{Lemma 2.9} $\{T_k\}$ is Cauchy in mass norm{\rm .}\endproclaim

\demo{Proof} Let $\varepsilon > 0$. Choose a compact subset
$L$   contained in $ \hat{\Gamma} \setminus E$
such that ${\cal H}^2(\hat{\Gamma} \setminus L) < \varepsilon/(2N_1)$
and such that each component of $L$ meets $K$.
Since $\hat{\Gamma} \setminus \Gamma$ consists of a finite
number of
components each of which meets $K$, we only need to  choose $L$,
using Lemma 2.8,
to meet each component  $C$ of $\hat{\Gamma} \setminus \Gamma$
in a sufficiently large
connected set which meets $K$ at some point of $ C$.

We claim that there exists $s_0$ such that $s>s_0$ implies
that $L$ is contained in the union of the components of $D_s$ which
meet $K$.  Indeed, set $\eta = \inf_L  \psi$. Choose $s_0$ so that
$\varepsilon_{s_0} < \eta$ and  $s > s_0$.  Let $x \in L$ and let
$\Omega$ be the component of $D_s$ which contains $x$. By the
choice of
$\eta$, $D_s \supseteq L $ and so $\Omega$ contains the component
of
$L$ through $x$. Hence $\Omega$ meets $K$, by the construction of
$L$.
This gives the claim.

Now suppose that $s_0 < j < k$ and consider $T_k-T_j$. We
claim that
  $\hbox{\rm mult}(\Sigma_k-\Sigma_j ,x) =0$
  if $x\in L$.
  Assume this for the moment. Hence
   $\hbox{\rm supp}(T_k-T_j) \subseteq  \hat{\Gamma} \setminus L$.
  Therefore, since
  $|\hbox{\rm mult}(\Sigma_k-\Sigma_j ,x)| \leq N_k+N_j\leq2N_1$
  by Lemma 2.7,  we get
  $${\bf M}(T_k-T_j) \leq  2N_1  {\cal H}^2(\hat{\Gamma} \setminus
L)
  < \varepsilon. $$
  This gives the lemma.

  It remains to verify the last claim. Let $x \in L$. Let $\Omega$ be
the
  component of $D_k$ containing $x$. Let $\Omega_1$ be the
component of
  $D_j \setminus \tau_j$ containing $x$.  Then $\Omega_1 \subseteq
\Omega$.
  Since $b(\Sigma_k - \Sigma_j) =\mu_k -\mu_j$ is disjoint from
  $\Omega_1$, $\Sigma_k - \Sigma_j$ has constant multiplicity
  in $\Omega_1$. By the  construction of $s_0$, $\Omega_1$ meets
  $K$ and so there is a path in $\Omega_1$ to  the point in $\Gamma
\setminus E$.
  Now  we can argue just as in the proof of Lemma 2.7 to conclude
that the
  multiplicity of $\Sigma_k - \Sigma_j$ is $0$ in $\Omega_1$. This
yields
  the claim. \enddemo
 
Let $T= \lim T_k$. By Lemma 2.9,  $T$ is the limit in mass norm
of the normal  currents $T_k$  and therefore is a
flat current in $\Bbb C^n$ (see [F]).
From $bT_k= [\sigma]+ [\mu_k]$ we conclude that
$bT= [\sigma]+ S$, where $S$ is a $1$-current supported
on $E$, since $\mu_k \subseteq Q_k$ and $\bigcap Q_k =E$.
Now, since  $T$ is flat,  $bT$ is flat. Hence $S$ is flat.
 Since $E$ is contained in the $1$-manifold $\Gamma$
 with   ${\cal H}^1(E)=0$ and since
 $dS=0$, we conclude, by Federer's support theorem, that $S=0$.
Hence $bT=[\sigma]$.

Now let $\phi$ be a holomorphic (1,0)-form in $\Bbb C^n$.
We have
$$\int_\sigma \phi= bT(\phi)=T(d\phi)=
\lim T_k(d\phi)= \lim \int_{\Sigma_k} d\phi.$$
But $d\phi=0$ on the $1$-variety $\hat{\Gamma} \setminus \Gamma
\supseteq \Sigma_k$.  Therefore
$\int_\sigma \phi=0$.
Since  $b[V]= [\Gamma]$ as currents,
we have
$\int_\Gamma \phi= \int_V d\phi=       0$, as $d\phi=0$
on $V_{\rm reg}$. Finally since $\gamma=\Gamma + \sigma$,
we have   $\int_\gamma\phi=0$. This is the moment
condition. This completes the proof of Theorem 3.\hfill\qed\enddemo

\demo{{R}emark} The proof shows that Cases 1 and   2a are not
possible.
Hence only Case 2b can occur and then
$[\gamma ]= [\Gamma] +[\sigma ]=b[V] +bT= b([V] +T)$.
That  is, $[\gamma ]$ is the boundary of the (1,1) current
$[V] +T$ supported in  $\hat{\gamma}$. This  gives a replacement
to the simpler condition that $\gamma \sim 0$ in $\hat{\gamma}$,
when
$\gamma$ is real analytic. \enddemo

\section{Proof of Theorem 1} 

\nonumproclaim{Lemma 3.1}  Let
$M$ be a $3$\/{\rm -}\/manifold in $\Bbb C^3$ satisfying the hypotheses
of Theorem {\rm 1.}  Let $\phi$ be a complex linear
function on $\Bbb C^3 $ and
let $H_\lambda$ be the affine complex hyperplane
$\{ z \in \Bbb C^3: \phi(z)=\lambda \}$ for $\lambda \in \Bbb C${\rm .}
Then{\rm ,} for almost all  $\lambda${\rm ,}
$H_\lambda  \cap   M$ is {\rm (}\/empty or\/{\rm )}  a smooth oriented
$1$\/{\rm -}\/cycle  $\gamma_\lambda $
{\rm (}\/i{\rm .}e{\rm .}\ a finite set of disjoint closed curves\/{\rm )}  that
bounds a   positive holomorphic $1$\/{\rm -}\/chain in
$H_\lambda \setminus \gamma_\lambda${\rm .}\endproclaim

 \demo{Proof} For almost all $\lambda$, by
 Sard's theorem, $M \cap H_\lambda$ is a smooth $1$-cycle,
call it $\gamma_\lambda$, and $\gamma_\lambda$ carries
an orientation as a slice
of $M$ by the map  $z  \mapsto \phi(z)$. We view $H_\lambda$
 as a copy of $\Bbb
C^2$.
We  claim that $\gamma_\lambda$ satisfies the  linking condition
hypothesis for Theorem 3. Let $A$ be an algebraic curve in
$\Bbb C^2  =H_\lambda$  that is  disjoint
from $\gamma_\lambda$. We can also view
$A$ as an algebraic curve in
$\Bbb C^3 \supseteq H_\lambda$ that is  disjoint from $M$.
By Lemma 1.1, ${\rm link}(A,M)$, with $A$   a curve in $\Bbb
C^3$,
agrees with
 ${\rm link}(A, \gamma_\lambda)$, with $A$   a curve in
$\Bbb C^2$.
 As  ${\rm link}(A,M) \geq 0$ by the hypothesis of Theorem 1,
 ${\rm link}(A,\gamma_\lambda) \geq 0$. Now Theorem  3
 implies that $\gamma_\lambda$  bounds
 a   positive holomorphic $1$-chain in
$H_\lambda \setminus \gamma_\lambda$. \enddemo

\demo{Proof of Theorem {\rm 1}}
We want to show that $M$  is
maximally complex. This means we need to show [HL] that
$\int_M \psi =0$ for every global $(p,q)$-form $\psi$ with
$p+q=3$ and $|p-q| >1$. This yields either $(3,0)$ or $(0,3)$.  By
complex conjugation, we can thus assume
that $\psi$ is a $(3,0)$-form; i.e.,
it suffices to show that
$$    \int_M \alpha   \, dz_1 \wedge   dz_2 \wedge dz_3
=0\leqno(3.1)$$
for all smooth functions $\alpha$ on $M$. Without loss of
generality  we can assume that $\alpha \in {\cal C}_0^\infty (\Bbb
C^3)$.
Using the inverse Fourier transform, we can write
$$ \alpha(z) = \int_{\Bbb C^3} e^{i(z\cdot \zeta+\overline{z\cdot
\zeta})}   \beta(\zeta) d{\cal L}^6(\zeta),$$
where
$z \cdot \zeta= z_1\zeta_1+z_2\zeta_2+z_3\zeta_3$,
${\cal L}^6$ is Lebesgue measure on $\Bbb C^3$, and
 $ \beta \in \cal S$\break ($\cal S$ is the class of ``rapidly decreasing''
functions).
Putting this integral expression for $\alpha$
into $(3.1)$  and interchanging the order of  integration,
we see that it suffices to show that
$$    \int_M   e^{i(z\cdot \zeta+\overline{z\cdot
\zeta})} dz_1 \wedge   dz_2 \wedge dz_3 =0\leqno(3.2)$$
for $\zeta \neq 0$, $\zeta \in \Bbb C^3$.

Now fix $\zeta \neq 0$. We introduce a
complex linear change of variable in $\Bbb
C^3 $  so that  $w_1= z \cdot (i\zeta), w_2, w_3$  are the new
variables.
Viewing  $w_1, w_2, w_3$ as functions on $\Bbb C^3$, we get
$dw_1 \wedge   dw_2 \wedge dw_3 =
c dz_1 \wedge   dz_2 \wedge dz_3 $
for some complex constant  $c \not =0$.
Thus $(3.2)$  is equivalent to the equation $I=0$ where
$$ I=  \int_M   e^{ w_1 - \bar{w_1}} dw_1 \wedge   dw_2 \wedge
dw_3 .$$
Let $\Sigma$ be a rectifiable $4$-chain in $\Bbb C^3$ such that
$b[\Sigma]=
[M]$---for example, a cone in $\Bbb C^3$ with base $M$.
Applying  Stokes' theorem to $I$ we get
$$I= \int_{\Sigma}
d(e^{ w_1-\bar{w_1}} dw_1 \wedge   dw_2 \wedge dw_3 )=
-\int_{ \Sigma}
  e^{ w_1-\bar{w_1}}
  d\overline{w_1} \wedge dw_1 \wedge   dw_2 \wedge dw_3 .$$
  Now slice $\Sigma$ by the map
  $w=(w_1,w_2,w_3) \mapsto w_1$ and get
  $$I= -\int_{\lambda \in \Bbb C} e^{\lambda -\bar{\lambda}}
  \left(\int_{\Sigma_\lambda} dw_2 \wedge dw_3\right) d\lambda
\wedge
  d\bar{\lambda},$$
  where $\Sigma_\lambda$  is the slice of $\Sigma$. Set
  $$ J_\lambda = \int_{\Sigma_\lambda} dw_2 \wedge dw_3 .$$
  It suffices to show that $J_\lambda=0$ for almost all $\lambda$.

  Let $H_\lambda$ be the affine hyperplane $\{w_1=\lambda\}$.
Then
  by Lemma 3.1, for almost all $\lambda$,
  $H_\lambda \cap M= \gamma_\lambda$ is a  $1$-cycle for which
there
  exists a holomorphic $1$-chain
   $[V_\lambda]$ with $b[V_\lambda]=[\gamma_\lambda]$.
   Also
  $b[\Sigma_\lambda]=[\gamma_\lambda]$ a.e.; this is because (a) by
a remark
  of  Harvey and Shiffman
  ([HS, 1.3.9, p.~567]), slicing commutes   a.e. with the $d$ operator, (b)
  $bT = -dT$ for $p$ currents with $p$ even, and (c) $\Sigma$ and
  $\Sigma_\lambda$ are a $4$-current  and a $2$-current,
respectively.
   Therefore,
  by two applications of Stokes' theorem, 
  $$J_\lambda= \int_{\Sigma_\lambda}dw_2 \wedge dw_3=
  \int_{\Sigma_\lambda}d(w_2 \wedge
dw_3)=\int_{\gamma_\lambda}
   w_2 \wedge dw_3=
  \int_{V_\lambda} dw_2 \wedge dw_3.$$
  Since  $[V_\lambda]$ is a $(1,1)$-current, the last integral equals
$0$.
  This completes the proof that $M$ is maximally complex. The
  following lemma then yields Theorem 1. \enddemo

  \demo{{R}emark} The idea of applying the Fourier transform in
connection
  with slicing is due to Globevnik and Stout [GS]. They  used it to show
  that if a function   satisfies the Morera property  on
  the boundary of a domain  in $\Bbb C^n$, then it
  satisfies the weak tangential
  Cauchy-Riemann equations.\enddemo

   \nonumproclaim{Lemma 3.2} Let $M$ be a {\rm (}\/not necessarily connected\/{\rm )}
   smooth compact oriented  $k$\/{\rm -}\/dimensional
   manifold  in $\Bbb C^n$  and suppose that
   $M$ is maximally complex{\rm .}
   Let
   $$T = \sum n_j [V_j]$$
   be the unique holomorphic $s$\/{\rm -}\/chain  {\rm (}\/with  $2s-1=k${\rm )}
   in $\Bbb C^n \setminus M$
   such that $bT=[ M]$ whose
   existence is given by the Harvey\/{\rm -}\/Lawson theorem{\rm .}
   If $M$ satisfies the linking condition{\rm ,} then
   $T$ is positive{\rm ;} i{\rm .}e{\rm .,}
   $n_j >0$
   for all $j${\rm .}
\endproclaim

   \demo{Proof}  We proceed by induction on odd $k$. The case $k=1$
   is   Lemma  2.1. Assume that $k \geq 3$.
   Let $F$ be a complex linear function that
   is not   locally constant on any of the varieties $V_j$.
    We consider   slices of
   rectifiable  currents by the map
   $F:\Bbb C^n \rightarrow \Bbb C$.
   Put $H_\lambda = \{z \in \Bbb C^n : F(z)=\lambda\}$, a
   complex hyperplane.
   For almost all $\lambda$
   we have that $\langle [V_j],F, \lambda\rangle  = [V _j \cap H_\lambda] $,
    $\langle M, F , \lambda\rangle  =[M \cap H_\lambda]$
   and $M \cap H_\lambda$ is a smooth oriented $(k-2)$-manifold
   satisfying the linking condition in $H_\lambda=\Bbb C^{n-1}$.

   Fix an index $m$. We can choose  $\lambda$ such that
   $V_m \cap H_\lambda$ contains an analytic branch $W$ (of
   complex dimension $s-1$) such that $W$ is not contained in any of
the
   $V_j$  for $ j \not = m$. We can choose $\lambda$ such that
   $\langle [V_j],F, \lambda\rangle  = [V_j  \cap H_\lambda] $ for all $j$,
   $\langle M, F , \lambda\rangle  =[M \cap H_\lambda]$,
   and $M \cap H_\lambda$ is a smooth oriented $(k-2)$-manifold
   satisfying the linking condition in $H_\lambda=\Bbb C^{n-1}$.
   Also $M \cap H_\lambda$ is maximally complex,
   since  $[M \cap H_\lambda]$ is the boundary of the
   holomorphic chain $S= \langle T,F,\lambda\rangle $ (see the next paragraph).
    Therefore, by induction,  the
    holomorphic $(s-1)$-chain   $S$ has {\it  positive\/}
   multiplicities.

   We have  $[M] =   bT$ and since slicing commutes almost
everywhere
   with  the boundary operator, we can also choose $\lambda$ so that
   $ \langle bT,F,\lambda\rangle =b\langle T,F,\lambda\rangle $. Hence
   $[M \cap H_\lambda]= b(\sum n_j [V_j \cap H_\lambda])$.
   And  so $S$ is the chain $\sum n_j [V_j \cap H_\lambda]$, where
   some of the  $V_j \cap H_\lambda$ may  be empty or reducible.
   Thus
   $n_m   [W]$ is one of the terms in the holomorphic $(s-1)$-chain
$S$
   when $S$ is written as a sum of irreducible varieties with integral
   multiplicities.
   Hence, $S$ being positive,  $n_m > 0$. \enddemo

   \demo{{R}emark} In order to proceed by induction,  we need
   to consider disconnected $M$.  However when $M$ is connected,
    one can argue more directly, since
   in that case , the Harvey-Lawson result for maximally complex
$M$
   is that
   $[M] = \pm b [V]$, for $V$ an irreducible $s$-variety in
   $\Bbb C^n \setminus M$.
   The linking condition then obviously implies
 that    $[M] =  b [V]$.\enddemo

\section{Proofs of Theorem 4 and 2}
\demo{Proof of Theorem {\rm 4}}
 It is clear that
$x \in  \hbox{\rm supp }T \setminus \gamma$ implies
$$ {\rm link}(\gamma, A)  > 0.$$
For the converse suppose that
$x \not \in \hbox{\rm supp }T \setminus \gamma$.
We shall show that
there exists an algebraic hypersurface  $A $ in $\Bbb C^n$ such that $x\in A$ and
$$ {\rm link}(\gamma,A)  = 0.$$
First suppose that  $x \not\in \hat{\gamma}$. Then there exists a
polynomial $P$ such that $P(x)=1 > ||P||_{\hat{\gamma}}$.
Then $A= {\bf Z}(P-1)$ gives the desired algebraic hypersurface.

Finally  suppose that $x \in \hat{\gamma} \setminus
\{\gamma \cup \hbox{\rm supp }T\}$.
Then, by Lemma 1.5, there exists a polynomial
$P$ in
   $\Bbb C^n$
   such that
    $P(x) \!=\!0$ and $P\! \neq\! 0$ on $\hbox{\rm supp}\,T\, \cup\, \gamma $.
   Hence $A= {\bf Z}(P)$
   gives an algebraic hypersurface such that $x \in A$ and
   ${\rm link}(\gamma,A)   =0$, since
   $P \neq 0$ on $\hbox{\rm supp }T \cup \gamma $.
   This gives Theorem 4.
\enddemo

   We  next show that Theorem 4 implies Theorem 2: It is clear that\break
$x \in \hbox{\rm supp }T$ implies
$$ {\rm link}(M ,A)  > 0.$$
For the converse take $x \not \in \hbox{\rm supp }T$
and choose an affine complex linear hyperplane
$H_\lambda=\{F=\lambda\}$ through $x$
such that $H_\lambda \cap M$ is an
oriented 1-manifold $\gamma$ so that
$[\gamma] = T_\lambda$, the $\lambda$-slice of $T$ by $F$.
We view $H_\lambda$ as a copy of $\Bbb C^{2}$
and apply Theorem 4 (taking $T_\lambda$ as the $T$ in Theorem 4).
Since  $x \not\in \hbox{\rm supp }T$,
we conclude that there exists an algebraic curve $A$ in $H_\lambda$
such that $x \in A$, $A \cap \gamma = \emptyset$ and
$$ {\rm link}(\gamma ,A)  = 0.$$
Then $A$ is an algebraic curve in $\Bbb C^3$ such that
$A \cap M= \emptyset$   and,  by Lemma 1.1,
$$ {\rm link}(M,A)  = {\rm link}(\gamma ,A) =0.$$
This yields Theorem 2.  \hfill\qed

\section{Proof of Theorem 5}

\nonumproclaim{Lemma 5.1}  Let $T$ be  an $\Bbb  R$\/{\rm -}\/linear subspace of
$\Bbb C^n$
of odd real dimension $k>3$ that is not maximally complex{\rm .}
Then there exists a complex  linear hyperplane $H$ {\rm (}\/through 0\/{\rm )} in
$\Bbb C^n$  such that $T  \cap H$ is not
maximally complex{\rm .}\endproclaim

\demo{Proof} Let $E = T \cap iT$ be the maximal complex
linear subspace of $T$.
Let $F$ be the Hermitian orthogonal complement of $E$. Then
 $E$ and $F$ are
complex linear subspaces of $\Bbb C^n$ such that
$\Bbb C^n= E \oplus  F$.  Hence
$$T= E \oplus  S$$
 where $S= T \cap F$ is totally real; i.e. $S \cap iS=\{0\}$.
 Since $T$ is not maximally complex
 and of odd dimension,
   $\dim S \geq 3$.

 We consider two cases:
\begin{itemize}
 \item[(a)] $E \not = \{0\}$,
 or
\item[(b)]  $E= \{0\}$, i.e.\ $T=S$ is totally real.
\end{itemize}

 {\it Case} (a) .  Take $u \not = 0$ in $E$ and set $H $
 equal to  the Hermitian
 orthogonal complement of $\Bbb C[u]$, the
 $\Bbb C$-span of $u$. Then $T \cap H= E' \oplus S$,
 where $E' = E \ominus u$ and
 $i(T \cap H)= E' \oplus iS$. Hence
 $(T \cap H)\ \cap  \ i(T \cap H)= E' $
   has real codimension in $T \cap H$ equal to $\dim S$,  which is
 greater than or equal to $3$,  and so $T \cap H$ is not
 maximally complex.

 {\it Case} (b). For any  complex  linear hyperplane
 $H$, $T \cap H$ is totally real, since $T$ is totally
 real. Since  $\dim T \cap H  \geq  k-2 \geq 3$ as $k \geq 5$,
 it follows that    $T \cap H$ is not
 maximally complex.

 \medbreak{\it Remark}. The condition that a linear space be totally real is an
``open''
 condition. This means that for all
  $\Bbb  R$-linear subspaces  $T'$    of $\Bbb C^n$ sufficiently close
to
   $T$ and with the same dimension and for all
   complex hyperplanes $H'$ sufficiently
 close to the hyperplane constructed in the proof of the lemma,
  $H' \cap T'$ is also not maximally complex. 

\demo{Proof of Theorem {\rm 5}} Assuming that $M$
satisfies the linking condition, we prove the maximal complexity
by induction on (odd) $k$, starting with $k=3$.
We have already done the case $k=3$ and $n=3$
as Theorem 1.
 Assume that $k=3$ and $n>3$ and fix $x \in M$.
  We claim that $T_x(M)$,
  the tangent space of $M$ at $x$, is maximally complex.
  There exists a complex linear map
   $\phi:\Bbb C^n \rightarrow \Bbb C^3$
   such that $x$ is
a regular point of $\phi|M$ and $\phi^{-1}(\phi(x)) \cap M =\{x\}$
(see Harvey [H, proof of Lemma 3.5, p.~349]).
Let $\cal M$ denote the image $\phi(M)$; $\cal M$ is a an oriented
immersed 3-manifold with
singularities (the ``scar set'') in $\Bbb C^3$ such that
integration over
$M$  gives a current $[{\cal M}]$ such that
$\phi_*([M] )= [{\cal M}]$
(see Harvey ([H, p.~368]).

 Let $H$ be a complex hyperplane in $\Bbb C^3$
 such that (a) $L= \phi^{-1}(H)$
is a complex hyperplane in $\Bbb C^n$
  that intersects  $M$ in a smooth  $1$-cycle
$\tilde{\gamma}$  and (b) $H \cap \cal M$ is
 a   $1$-cycle $\gamma =\phi(\tilde{\gamma})$.
 As $M$ satisfies the linking condition,
our previous arguments show that $\tilde{\gamma}$ satisfies the
linking
condition in $L$ (with $L$ viewed as a copy of $\Bbb C^{n-1}$).  Hence
$\tilde{\gamma}$ bounds a holomorphic $1$-chain $\tilde{V}$.
It follows that $\gamma$ bounds the holomorphic $1$-chain
$V= \phi(\tilde{V})$.
As (a) and (b) hold for almost all hyperplanes $H$,
 we conclude, as   in the proof of Theorem 1, that
$$ \int_{  {\large {\cal M} } }  \alpha \  dz_1 \wedge dz_2 \wedge
dz_3= 0$$
for all  ${\cal C}^\infty$  functions $\alpha$  on $\Bbb C^3$.  This
implies that
$T_y(\cal M)$ is maximally complex  at points $y \in \cal M$
where $\cal M$ is smooth,
since $\alpha$ can be chosen  to have support in
arbitrarily small neighborhoods of $y$.
It follows that $T_x(M)$ is maximally complex, because $x$ is a
regular point of $\phi|M$. As $x$ is an arbitrary point of $M$,
we conclude
that $M$ is maximally complex. By Lemma 3.2 this completes the
 case $k=3$.

Now consider the case $k >3$.  We argue by contradiction and
suppose that  $M$ is not maximally complex. It follows that
there exists  a point $p \in M$ such that
$T_p(M)$, the tangent space of $M$ viewed as a real linear subspace
of
$\Bbb C^n$, is not  maximally complex.
Then by  Lemma 5.1, there is a complex hyperplane $H$ in $\Bbb
C^n$
such that $T_p(M) \cap H_p$ is not maximally complex.
Suppose, for the moment,  that the translate
$H_p=H+p$  of $H$ through $p$ in $\Bbb C^n$, given say  by
$H_p ={\bf Z}(F)$ where $F$ is an affine complex linear function,
is such  that $Q =M  \cap H_p$ is a smooth $(k-2)$-manifold,
   oriented as the slice of the map  $F:M \rightarrow \Bbb C$.
Then $T_p(Q) =T_p(M) \cap H   $ is not maximally complex.
On the other hand,
 the linking condition for $M$ implies that $Q$ satisfies the
linking condition in $H_p$ (viewed as a copy of $\Bbb C^{n-1}$).
Hence, by
induction on $k$, since $k-2 \geq 3$, we conclude that $Q$ is  a
maximally
complex manifold. Therefore $T_p(Q) $ is  a maximally complex
linear space.
This is a contradiction.

 From this we can conclude that $M$ is maximally
  complex---except for the assumption above that
$Q =M  \cap H_p$ is a smooth $(k-2)$-manifold. If this is not the
case, we
can, by Sard's theorem, translate $H_p $ a small amount so that
the intersection with $M$ becomes smooth. Our  above argument,
together with the remark
after  Lemma 5.1, then show  that
$T_q(M) $ is maximally complex for some
 $q \in M$ arbitrarily close to $p$. It follows that, in the limit,
    $T_p(M) $ is maximally complex. Again Lemma 3.2 yields the first
part
    of Theorem 5 for general $k>0$.

    Finally,  we  prove by induction on $k$ that if
    $x \in \Bbb C^n \setminus M$  and  $x \not \in  \hbox{\rm supp }T$,
    then there exists an algebraic subvariety $A$
    of $\Bbb C^n$ with $x \in A$  and
    $${\rm link}(M,A)  = 0.$$
 The case $k=3$ is just Theorem 2.

Suppose that $k >3$. We argue as in the proof of Theorem 2.
Choose an affine complex linear hyperplane
$H_\lambda=\{F=\lambda\}$ through $x$
such that $H_\lambda \cap M$ is an
oriented $(k-2)$-manifold $Q$ so that
$[Q] = b(T_\lambda)$, where $T_\lambda$ is the slice of $T$
by $F$.  We view $H_\lambda$ as a copy of $\Bbb C^{n-1}$.
As $x \not \in Q$,
we can apply the induction hypothesis to get an algebraic
subvariety $A$ of $H_\lambda$ such that $A $ is disjoint from  $ Q$,
$x \in A$ and $ {\rm link}(Q ,A)  = 0$, where the linking number is
computed in $\Bbb C^{n-1}$.
Then $A$ is an algebraic subvariety in $\Bbb C^n$ such that
$A$ is disjoint from $M$ and
$ {\rm link}(M,A)  = 0$. This  gives the theorem. \hfill\qed

\AuthorRefNames [AW]

 \end{document}